\documentclass{article}
\usepackage{amsmath,amssymb,bbm,graphics,chicago,graphicx,latexsym}
\title{Likelihood ratio type two-sample tests for current status data}
\author{Piet Groeneboom\\
{\it Delft University of Technology}}

\def\d{\delta}

\def\th{\theta}

\def\P{{\mathbb P}}
\def\G{{\mathbb G}}
\def\H{{\mathbb H}}

\def\a{\alpha}
\def\b{\beta}
\def\R{\mathbb R}
\def\dd{\Delta}

\def\d{\delta}

\def\G{{\mathbb G}}
\def\H{{\mathbb H}}

\def\P{{\mathbb P}}

\def\l{\lambda}
\def\labda1{\lambda_1}
\def\labda2{\lambda_2}

\def\e{\varepsilon}
\def\f{\phi}
\def\t{\tau}

\def\s{\sigma}

\def\comment#1{\relax}

\def\=in{\mathop{\rm =}}

\def\eop{\hfill\mbox{$\Box$}\newline}

\newtheorem{theorem}{Theorem}[section]

\newtheorem{lemma}{Lemma}[section]

\newtheorem{remark}{Remark}[section]

\begin{document}
\maketitle

\noindent ABSTRACT. We introduce fully nonparametric two-sample tests for testing the null hypothesis that the  samples come from the same distribution if the values are only indirectly given via current status censoring. The tests are based on the likelihood ratio principle and allow the observation distributions to be different for the two samples, in contrast with earlier proposals for this situation. A bootstrap method is given for determining critical values and asymptotic theory is developed. A simulation study, using Weibull distributions, is presented to compare the power behavior of the tests with the power of other nonparametric tests in this situation.

\vspace{5mm}

\section{Introduction}
\label{intro}
\setcounter{equation}{0}
At the beginning of the vast amount of research on right-censored data, there was much interest in two-sample tests for right-censored data, like the Gehan test, log rank test, Efron's test, etc. For example,  \shortciteN{gehan:65} considers the testing problem of testing $F_1\equiv F_2$ against the alternative $F_1<F_2$, and gives a permutation test for this testing problem.

Permutation tests for the two-sample problem with interval censored data have been considered in 
\shortciteN{peto+peto:72}. Since they rely on the permutation distribution,
such tests can only be used when the censoring mechanism is the same in both
samples. One of the referees of this paper asked the interesting question whether permutation tests of this type, considered as conditional tests, might be asymptotically independent of the observation distributions in the two samples, in analogy with results in \shortciteN{neuhaus:93} for two-sample tests in the presence of right censoring.
I do not know the answer to this question (current status censoring is very different from right censoring!), but preliminary results indicate that this method gives very variable estimates of the critical values for moderate sample sizes and therefore cannot be used for these sample sizes. The bootstrap method we propose for computing the critical values does not suffer from this drawback, see section \ref{section:sim_study}.

The maximum likelihood estimator for interval censored data is considered in more detail in \shortciteN{peto:73}, where it is suggested that pointwise standard errors for the survival curve
can be estimated from the inverse of the Fisher information. However, we know by now for a long time that this is not correct if we sample from continuous distributions; the pointwise asymptotic distribution is not normal, and the asymptotic variance is not given by the the inverse of the Fisher information, see, e.g., \shortciteN{GrWe:92} (I owe this observation on \shortciteN{peto:73} to Peter Sasieni).

Other tests have been considered in, e.g., \shortciteN{andersen_ronn:95} and \shortciteN{sun:06}, where also references to earlier work by the latter author can be found. They are based on certain functionals of the distributions which will be different from zero for some alternatives (mostly of the type of ``shift alternatives"). Similar tests have been considered in \shortciteN{zhang:01} and \shortciteN{zhang:06} for panel count data, where pseudo maximum likelihood estimators are used. Specialized to our present problem, this leads to tests of the same type as the tests in \shortciteN{andersen_ronn:95} and \shortciteN{sun:06}.

We consider here rather different types of tests which are likelihood ratio based tests for testing that two samples come from the same distribution, if current status censoring is present. A test of this type is considered in Chapter 3 of \shortciteN{vlad:03}, where the null hypothesis of equality of the distribution functions $F_1$ and $F_2$, generating the first and second sample, respectively, is tested against Lehmann alternatives of the form
\begin{equation}
\label{Lehmann_alt}
F_2(t)=F_1(t)^{1+\th},\,\th\in(-1,\infty)\setminus\{0\}.
\end{equation}

Here we prefer to test the null hypothesis of equality of $F_1$ and $F_2$ just against the more general alternative that they are not equal. Note that in testing against the Lehmann alternatives (\ref{Lehmann_alt}), we have to estimate $F_1$ and $\th$, whereas in the more general testing problem we have to estimate both $F_1$ and $F_2$ nonparametrically.

We will assume the usual conditions for the current status model with continuous distributions, as stated on p.\ 35 of \shortciteN{GrWe:92}: $(X_1,T_1),\dots,(X_m,T_m)$ and $(X_{m+1},T_{m+1}),\dots,(X_{N},T_{N})$, $N=m+n$, are two independent samples of random variables in $\R^2$, where $X_i$ and $T_i$ are independent, with, respectively, continuous distribution functions $F_1$ and $G_1$ in the first sample and continuous distribution functions $F_2$ and $G_2$ in the second sample. We call the $X_i$ the ``hidden" variables and the $T_i$ the observation variables. Note that we allow the distribution functions $G_1$ and $G_2$ of the observation variables to be different in the two samples. In the current status model, the only observations which are available to us are the pairs
$$
(T_i,\dd_i),\qquad\dd_i=1_{\{X_i\le T_i\}},
$$
so we do not observe $X_i$ itself, but only its ``current status" $\dd_i$. In this situation, we want to test the null hypothesis that the distribution functions of the hidden variables are the same in the two samples.

We first discuss what a simple likelihood ratio test would look like. Under the null hypothesis we have to maximize
$$
\sum_{i=1}^N \left\{\dd_i\log F(T_i)+(1-\dd_i)\log\left(1- F(T_i)\right)\right\},\qquad N=m+n,
$$
over all distribution functions $F$, and without the restriction of the null hypothesis we have to maximize
\begin{align*}
&\sum_{i=1}^m \left\{\dd_i\log F_1(T_i)+(1-\dd_i)\log\left(1- F_1(T_i)\right)\right\}\\
&\qquad\qquad+\sum_{i=m+1}^N \left\{\dd_i\log F_2(T_i)+(1-\dd_i)\log\left(1- F_2(T_i)\right)\right\}
\end{align*}
over all pairs of distribution functions $(F_1,F_2)$.

This means that under the null hypothesis the MLE (maximum likelihood estimator) is given by the left-continuous slope of the greatest convex minorant of the cusum diagram of the points $(0,0)$ and the points
\begin{equation}
\label{cusum_combined}
\left(i,\sum_{j\le i}\dd_{(j)}\right),\,i=1,\dots,N.
\end{equation}
using a notation, introduced in \shortciteN{GrWe:92}. Here $\dd_{(j)}$ denotes the indicator corresponding to the $j$th order statistic $T_{(j)}$. Without the restriction of the null hypothesis the MLE of $F_1$ is given by the left-continuous slope of the greatest convex minorant of the cusum diagram of the points $(0,0)$ and the points
\begin{equation}
\label{cusum_first}
\left(i,\sum_{j\le i}\dd_{(j1)}\right),\,i=1,\dots,m,
\end{equation}
where $\dd_{(j1)}$ is the indicator corresponding to $j$th
order statistic $T_{(j1)}$ of the first sample. Similarly the MLE of $F_2$ is given by the left-continuous slope of the greatest convex minorant of the cusum diagram of the points $(0,0)$ and the points
\begin{equation}
\label{cusum_second}
\left(i,\sum_{j\le i}\dd_{(j2)}\right),\,i=1,\dots,n,
\end{equation}
where $\dd_{(j2)}$ is the indicator corresponding to $j$th
order statistic $T_{(j2)}$ of the second sample.

Let the MLE of $F_1$ ($=F_2$) under the null hypothesis be given by $\hat F_N$, and let the MLE of the pair $(F_1,F_2)$ without the restriction of the null hypothesis be given by
$$
\left({\hat F}_{N1},{\hat F}_{N2}\right).
$$
Then the log likelihood ratio test statistic is given by:
\begin{align}
\label{LR_statistic}
&\sum_{i=1}^m \left\{\dd_i\log \frac{{\hat F}_{N1}(T_i)}{{\hat F}_N(T_i)}+(1-\dd_i)\log\frac{1- {\hat F}_{N1}(T_i)}{1- {\hat F}_N(T_i)}\right\}\nonumber\\
&\qquad\qquad+\sum_{i=m+1}^N \left\{\dd_i\log \frac{{\hat F}_{N2}(T_i)}{{\hat F}_N(T_i)}+(1-\dd_i)\log\frac{1- {\hat F}_{N2}(T_i)}{1- {\hat F}_N(T_i)}\right\},
\end{align}
where the terms with coefficients $\dd_i$ and $1-\dd_i$ are defined to be zero if $\dd_i$ and $1-\dd_i$ are zero, respectively.

Although we take this statistic as our inspiration, we first study a statistic somewhat similar to this LR statistic, based on maximum smoothed likelihood estimators (MSLEs), introduced in \shortciteN{piet_geurt_birgit:10}. One of the reasons is that the asymptotic analysis of the original LR statistic is rather involved; the difficulty in analyzing the limit properties of (\ref{LR_statistic}) lies in the problem of finding a normalization making it an asymptotic pivot under the null hypothesis. One also has to deal with the non-standard asymptotics, which derives from the fact that the statistic is based on (non-linear) isotonic estimators which satisfy an order restriction. These non-standard features also turn up in the limit behavior. Another reason is that the MSLE leads to more powerful tests for models, commonly used in this type of comparisons. This will be illustrated by a simulation study for a two parameter Weibull distribution, also used in \shortciteN{andersen_ronn:95} in a simulation study to check the power of their proposed test.

Maximum smoothed likelihood estimators for current status data were studied in \shortciteN{piet_geurt_birgit:10}, where it was shown that, under some regularity conditions, the local limit distribution is normal (in contrast with the limit behavior of the original MLE). These estimators are obtained by first smoothing the observation distribution, for example by kernel estimators, and next maximizing the smoothed likelihood w.r.t.\ the distribution of the hidden variables. In this way the MSLE inherits smoothness properties of the estimate of the observation distribution and converges at a faster rate than the ``raw" MLE, which locally converges at rate $n^{-1/3}$ under the usual smoothness conditions on the underlying distributions. Further results on the MSLE can be found in \shortciteN{piet_geurt_birgit:10}.

\begin{figure}[!ht]
\begin{center}
\includegraphics[scale=0.4]{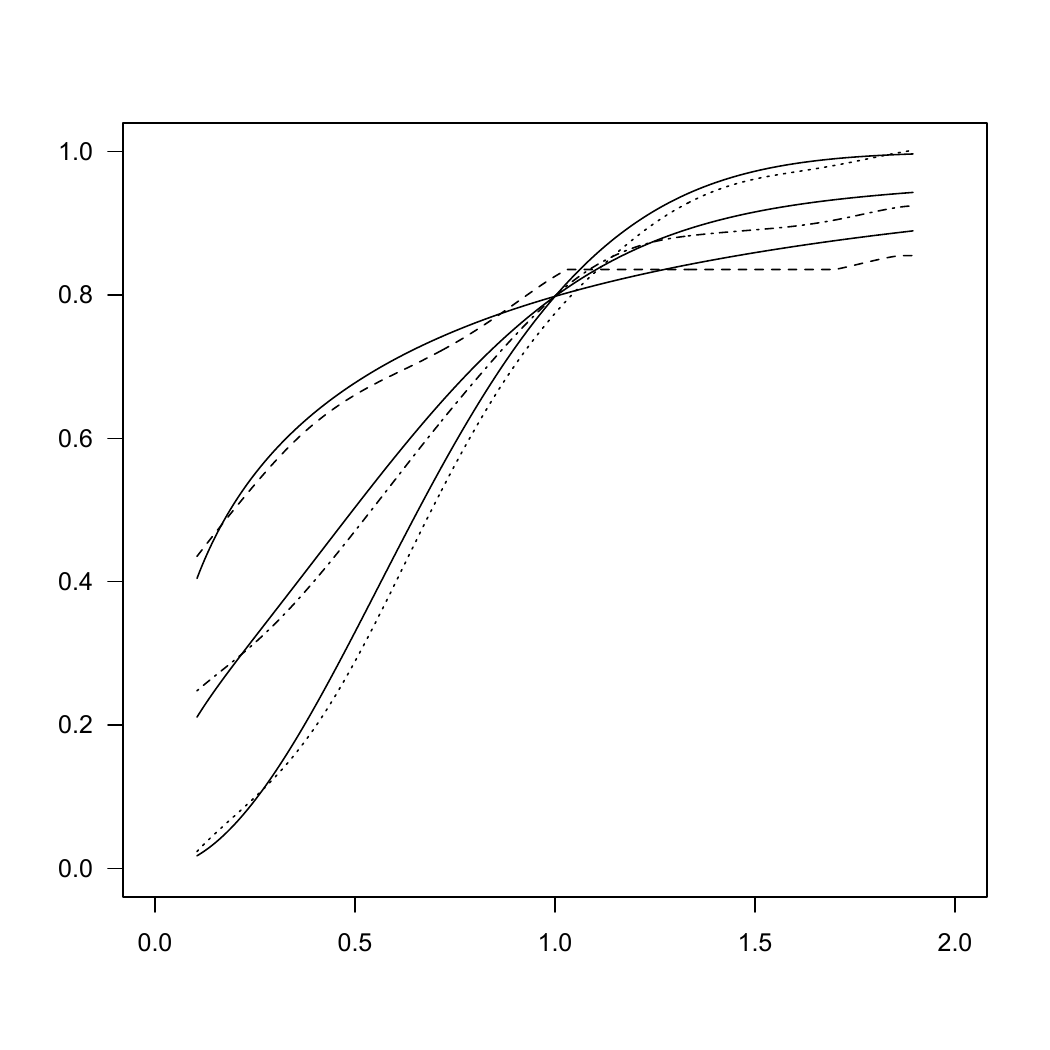}
\includegraphics[scale=0.4]{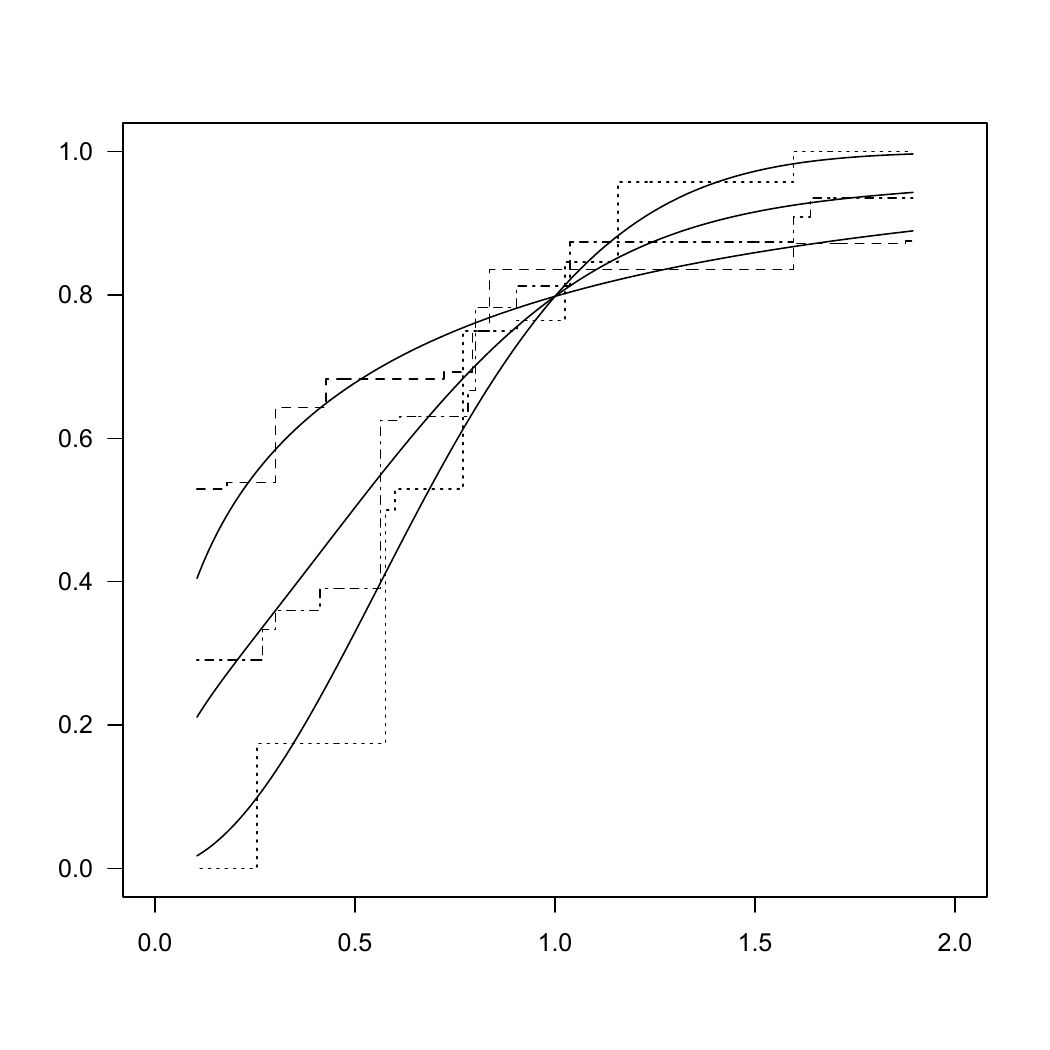}
\end{center}
\caption{MSLEs and MLEs on $[a,b]$ for samples of size $m=n=250$ from the Weibull densities (\ref{Weibull}). $G_1$ and $G_2$ are uniform on $[0,2]$, and the interval $[a,b]=[0.1,1.9]$. The left panel gives the MSLE estimates and the right panel the MLEs, where the dashed curves give the estimates for the first sample ($\a_1=0.5$), the dotted curves the estimates for the second sample ($\a_2=2$), and the dashed-dotted curves the estimates for the combined samples. The solid curves give the corresponding actual distribution functions for these three situations. The bandwidth for the computation of the MSLEs was $b_N=2N^{-1/5}\approx0.57708$, where $N=m+n=500$.}
\label{fig:MSLE+MLE_Weibull}
\end{figure}

A picture of the MSLE estimators and the MLE estimators for samples of size $250$ from two different Weibull distributions with densities
\begin{equation}
\label{Weibull}
\a_1\l x^{\a_1-1} e^{-\l x^{\a_1}},\qquad \a_2\l x^{\a_2-1} e^{-\l x^{\a_2}},\qquad\,x>0,\qquad\,\a_1=0.5,\,\a_2=2,\,\l=1.6,
\end{equation} 
respectively, where $\a_1=0.5$ holds for the first sample and $\a_2=2$ for the second sample, is shown in Figure \ref{fig:MSLE+MLE_Weibull}. 

\section{A likelihood ratio test, based on maximum smoothed likelihood estimators}
\label{section:MSLE}
\setcounter{equation}{0}
In order to avoid problems at the boundary, we restrict the domain on which we compute our test statistic to an interval $[a,b]\subset (0,M)$, where $[0,M]$ is assumed to be the support of the underlying densities, corresponding to the distribution functions $F_1$ and $F_2$ of the hidden variables. We consider the statistic $V_N$, similar to (\ref{LR_statistic}), and defined by
\begin{align}
\label{kernel_LR_statistic}
&V_N=\frac{2m}{N}\int_{t\in[a,b]} \left\{\tilde h_{N1}(t)\log \frac{{\tilde F}_{N1}(t)}{{\tilde F}_N(t)}+\bigl\{\tilde g_{N1}(t)-\tilde h_{N1}(t)\bigr\}\log\frac{1- {\tilde F}_{N1}(t)}{1- {\tilde F}_N(t)}\right\}\,dt\nonumber\\
&\qquad\qquad+\frac{2n}{N}\int_{t\in[a,b]} \left\{\tilde h_{N2}(t)\log \frac{{\tilde F}_{N2}(t)}{\tilde F_N(t)}+\bigl\{\tilde g_{N2}(t)-\tilde h_{N2}(t)\bigr\}\log\frac{1- {\tilde F}_{N2}(T_i)}{1- {\tilde F}_N(T_i)}\right\}\,dt
\end{align}
where ${\tilde F}_{N1}$, ${\tilde F}_{N2}$ and ${\tilde F}_N$ are the maximum smoothed likelihood estimators (MSLEs) for the first, second and combined sample, respectively, and $\tilde g_{Ni}$ and $\tilde h_{Ni}$ are kernel estimates of the relevant observation densities, defined below. As explained in \shortciteN{piet_geurt_birgit:10}, where the same type of MSLE for the current status model is defined, the MSLEs for the combined samples and the first and second sample are computed by replacing the cusum diagrams (\ref{cusum_combined}), (\ref{cusum_first}) and (\ref{cusum_second}) by the continuous cusum diagrams
\begin{equation}
\label{cusum_combined_smooth}
\left(\tilde G_N(t),\tilde H_N(t)\right),\,t\in\left[0,M\right], 
\end{equation}
\begin{equation}
\label{cusum_first_smooth}
\left(\tilde G_{N1}(t),\tilde H_{N1}(t)\right),\,t\in\left[0,M\right],
\end{equation}
and
\begin{equation}
\label{cusum_second_smooth}
\left(\tilde G_{N2}(t),\tilde H_{N2}(t)\right),\,t\in\left[0,M\right], 
\end{equation}
respectively, where $\tilde G_N$, $\tilde G_{Ni}$, $\tilde H_N$, $\tilde H_{Ni}$ and their derivatives are defined in the following way.

We first define the densities $\tilde g_{Ni}$ and $\tilde h_{Ni}$ on $[b_N,M-b_N]$ by
\begin{equation}
\label{smooth_dens}
\tilde g_{Ni}(t)=\int K_{b_N}(t-u)\,d\G_{Ni}(u),\qquad  \tilde h_{Ni}(t)=\int K_{b_N}(t-u)\,\d \,d\P_{Ni}(u,\d),
\end{equation}
Here $\G_{N1}$ is the empirical distribution function of the observations $T_1,\dots,T_m$ of the first sample and $\P_{N1}$ is the empirical distribution of the observations $(T_1,\dd_1),\dots,(T_m,\dd_m)$ of the first sample, with the analogous definitions of $\G_{N2}$ and $\P_{N2}$ for the second sample. The densities $\tilde g_N$ and $\tilde h_N$ are defined on $[b_N,M-b_N]$ by
$$
\tilde g_N=\a_N\tilde g_{N1}+\b_N\tilde g_{N2},\qquad \tilde h_N=\a_N\tilde h_{N1}+\b_N\tilde h_{N2},
\qquad\a_N=\frac{m}{N},\qquad\b_N=1-\a_N.
$$
The kernel $K_b$ is defined in the usual way by
$$
K_b(u)=\frac1b K(u/b),
$$
for a bandwidth $b>0$, where $K$ is a symmetric positive kernel with compact support. We consider symmetric positive polynomial-type kernels $K$, with compact support. In our simulation study we took
\begin{equation}
\label{triweight}
K(u)=\tfrac{35}{32}\left(1-u^2\right)^31_{[-1,1]}(u),
\end{equation}
the so-called triweight kernel.

For $t\in[0,b_N]$ and $t\in[M-b_N,M]$ we use a boundary kernel, defined by a linear combination of $K(u)$ and $uK(u)$. Other ways of bias correction at the boundary are also possible, but it seems necessary to use such a correction in order to obtain a reasonable behavior at the boundary. Using boundary kernels, we lose the simple property that the distribution function can be obtained by just integrating the kernel, and indeed the estimates of the distribution functions were obtained by numerically integrating the estimates of the densities (and not by integrating the kernels).
So we define
$$
\tilde G_{Ni}(t)=\int_0^t \tilde g_{Ni}(u)\,du,\qquad \tilde H_{Ni}(t)=\int_0^t \tilde h_{Ni}(u)\,du,
$$
$$
\tilde G_N=\a_N\tilde G_{N1}+\b_N\tilde G_{N2},\qquad \tilde H_N=\a_N\tilde H_{N1}+\b_N\tilde H_{N2},
$$
and use the corresponding numerical integrals in the continuous cusum diagrams (\ref{cusum_combined_smooth}) to (\ref{cusum_second_smooth}).

Note that the cusum diagrams (\ref{cusum_combined_smooth}) to (\ref{cusum_second_smooth}) are continuous analogues of the cusum diagrams (\ref{cusum_combined}) to (\ref{cusum_second}), since, for example, the left-continuous slope of (\ref{cusum_combined}) is the same as the left-continuous slope of the cusum diagram consisting of the set of points
$$
\left\{\left(\G_N(t),\H_N(t)\right),\,t\ge0\right\},
$$
where $\G_N$ is the empirical distribution function of the points $T_i,\,i=1,\dots,N$, and $\H_N$ is the empirical sub-distribution function of the points $T_i,\,i=1,\dots,N$, with $\dd_i=1$. However, the slopes of the greatest convex minorants of the continuous cusum diagrams (\ref{cusum_combined_smooth}) to (\ref{cusum_second_smooth}) are continuous functions of $t$ in contrast with the left-continuous slopes of the cusum diagrams (\ref{cusum_combined}) to (\ref{cusum_second}).

The following result shows that the test statistic $V_N$ is, for a suitable choice of the bandwidth, an asymptotic pivot under the null hypothesis of equality of the two distribution function $F_1$ and $F_2$ of the hidden variables in the two samples.

\begin{theorem}
\label{th:local_asymp_curstat_MSLE}
Let the test statistic $V_N$ be defined by (\ref{kernel_LR_statistic}), using a bandwidth $b_N$ such that $b_N\asymp n^{-\a}$, where $2/9<\a<1/3$. Furthermore, let $F$ stay away from zero and one on $[a,b]$ and have a bounded continuous second derivative $f'$ on an interval $(a',b')$ containing $[a,b]$, and let $g_1$ and $g_2$ be continuous densities which stay away from zero on $[a,b]$, with continuous bounded second derivatives on the interval $(a',b')$.
Let the log likelihood ratio statistic $V_N$, based on the MSLEs, be defined by (\ref{kernel_LR_statistic}). Then we have in probability, if the distribution functions of the hidden variables in the two sample are both equal to $F$ and $m/N\to\a\in(0,1)$, as $N\to\infty$,
\begin{equation}
\label{inprob_convergence}
N\sqrt{\frac{b_N}{b-a}}\left\{\bigl(V_N|T_1,\dots,T_N\bigr)-\frac{b-a}{Nb_N}\int K(u)^2\,du\right\}\stackrel{{\cal D}}\longrightarrow N(0,\s_K^2),
\end{equation}
where $N(0,\s_K^2)$ denotes a normal distribution with mean zero and variance
$$
\s_K^2=2\int\left\{\int K(u+v)K(u)\,du\right\}^2\,dv.
$$
\end{theorem}

\begin{remark}
{\rm
To say that (\ref{inprob_convergence}) holds in probability means that
$$
\P\left\{N\sqrt{\frac{b_N}{b-a}}\left\{V_N-\frac{b-a}{Nb_N}\int K(u)^2\,du\right\}\le x\Biggm|T_1,\dots,T_N\right\}\stackrel{p}\longrightarrow\Phi(x),
$$
for each $x\in\R$, where $\Phi$ is the standard normal distribution function and $\stackrel{p}\longrightarrow$ denotes convergence in probability.
}
\end{remark}

\begin{remark}
{\rm The restriction of the bandwidth to the range $N^{-1/3}\ll b_N\ll N^{-2/9}$ has the following motivation. The condition $b_N\gg N^{-1/3}$ is necessary for having the asymptotic equivalence of the MSLEs to ratios of kernel estimators (see Corollary 3.4 in \shortciteN{piet_geurt_birgit:10}), and $b_N\ll N^{-2/9}$ prevents the bias to enter, which causes  the asymptotic distribution of $V_N$ to become dependent on the observation densities $g_1$ and $g_2$. The bias term drops out if the observation densities $g_1$ and $g_2$ are the same in the two samples.

Nevertheless, we prefer to work with a larger bandwidth, at the cost of introducing a bias term, depending on the underlying distributions, as shown in Theorem \ref{th:SLR_test}. It turns out that this bias term does not bother us, if we compute the critical values by a bootstrap procedure, to be discussed in section \ref{section:crit_value_bootstrap}. The key to this is that the bias term is estimated automatically in the bootstrap resampling from a smooth estimate of $F$ and that the difference between this estimate of the bias and the bias is sufficiently small, as shown in the proof of Theorem \ref{th:bootstrap1}, so that we can replace it by the deterministic bias in the central limit theorem for the bootstrap test statistic.
}
\end{remark}

\begin{theorem}
\label{th:SLR_test}
Let the test statistic $V_N$ be defined by (\ref{kernel_LR_statistic}), using a bandwidth $b_N$ such that $b_N\asymp n^{-\a}$, where $1/5<\a\le 2/9$. Furthermore, let $F$ stay away from zero and one on $[a,b]$ and have a bounded continuous second derivative $f'$ on an interval $(a',b')$ containing $[a,b]$, and let $g_1$ and $g_2$ be continuous densities which stay away from zero on $[a,b]$, with continuous bounded second derivatives on the interval $(a',b')$.
Let the log likelihood ratio statistic $V_N$, based on the MSLEs, be defined by (\ref{kernel_LR_statistic}). Then we have in probability, if the distribution functions of the hidden variables in the two sample are both equal to $F$ and $m/N\to\a\in(0,1)$, as $N\to\infty$,
\begin{align*}
&N\sqrt{\frac{b_N}{b-a}}\left\{\bigl(V_N|T_1,\dots,T_N\bigr)-\frac{b-a}{Nb_N}\int K(u)^2\,du\right.\\
&\qquad\qquad\qquad\qquad\left.-\a_N\b_N\int_{t=a}^b\frac{f(t)^2\left\{g_1'(t)g_2(t)-g_2'(t)g_1(t)\right\}^2}{F(t)\{1-F(t)\}\bar g_N(t)g_1(t)g_2(t)}\,dt\left\{\int u^2 K(u)\,du\right\}^2b_N^4\right\}\\
&\stackrel{{\cal D}}\longrightarrow N(0,\s_K^2),
\end{align*}
where $\bar g_N$ is defined by:
$$
\bar g_N(t)=\a_Ng_1(t)+\b_Ng_2(t).
$$
and $N(0,\s_K^2)$ denotes a normal distribution with mean zero and variance $\s_K^2$ defined as in Theorem \ref{th:local_asymp_curstat_MSLE}.
\end{theorem}

\begin{remark}
{\rm If $b_N\asymp N^{-1/5}$ the situation becomes even more complicated. If the observation densities $g_1$ and $g_2$ are the same, we still get the asymptotic normality result, as shown in the following theorem. But if the densities $g_1$ and $g_2$ are different, extra non-negligible random terms enter because of the presence of the bias term. We will not discuss this further in the present paper.
}
\end{remark}

\begin{theorem}
\label{th:SLR_test_same_g}
Let the test statistic $V_N$ be defined by (\ref{kernel_LR_statistic}), using a bandwidth $b_N$ such that $b_N\asymp n^{-\a}$, where $1/5\le \a<1/3$. Furthermore,let $F$ stay away from zero and one on $[a,b]$ and have a bounded continuous second derivative $f^{\prime}$ on an interval $(a',b')$ containing $[a,b]$, and let $g_1=g_2$ be a continuous density which stays away from zero on $[a,b]$, with a continuous bounded second derivative on the interval $(a',b')$. Then we have in probability, if the distribution functions of the hidden variables in the two sample are both equal to $F$ and $m/N\to\a\in(0,1)$, as $N\to\infty$,
\begin{align*}
&N\sqrt{\frac{b_N}{b-a}}\left\{\bigl(V_N|T_1,\dots,T_N\bigr)-\frac{b-a}{Nb_N}\int K(u)^2\,du\right\}
\stackrel{{\cal D}}\longrightarrow N(0,\s_K^2),
\end{align*}
where $N(0,\s_K^2)$ denotes a normal distribution with mean zero and variance $\s_K^2$ defined as in Theorem \ref{th:local_asymp_curstat_MSLE}.
\end{theorem}

\begin{remark}
{\rm
We used a conditional formulation, since we will use conditional tests in our bootstrap approach, but the convergence in distribution will also hold in Theorems \ref{th:local_asymp_curstat_MSLE} to \ref{th:SLR_test_same_g}, if we do not condition on $T_1,\dots,T_N$.
}
\end{remark}

\section{The original LR test}
\label{section:LR}
\setcounter{equation}{0}
We return to the original LR test, using the MLEs, and confine ourselves to a heuristic discussion, since a complete treatment is still out of our grasp. As in the proof of Theorem \ref{th:local_asymp_curstat_MSLE}, we have:
\begin{align*}
&2\int_{[a,b]}\left\{\hat F_{N1}(t)\log\frac{\hat F_{N1}(t)}{\hat F_N(t)}+\bigl\{1-\hat F_{N1}(t)\bigr\}\log\frac{1-\hat F_{N1}(t)}{1-\hat F_N(t)}\right\}\,d\G_{N1}(t)\\
&\sim \int_{[a,b]}\frac{\bigl\{\hat F_{N}(t)-\hat F_{N1}(t)\bigr\}^2}{F(t)\bigl\{1-F(t)\bigr\}}\,dG_1(t),
\end{align*}
with a similar relation for the terms involving $\hat F_{N2}$. This motivates the study of integrals of the following type:
$$
E \int_a^b\frac{\bigl\{\hat F_N(x)-F(x)\bigr\}^2}{F(x)\{1-F(x)\}}\,dG(x).
$$

The local limit of the MLE of the combined samples under the null hypothesis, when the observation times $T_i$ in both samples is given by $G$ is given in the following theorem, given on p.\ 89 of \shortciteN{GrWe:92}.

\begin{theorem}
\label{th:cur_stat_convergence}
Let $t_0$ be such that $0<F(t_0),G(t_0)<1$, and let
$F$ and $G$ be differentiable at
$t_0$, with strictly positive derivatives $f(t_0)$ and $g(t_0)$, respectively.
Furthermore, let $\hat F_N$ be the MLE of $F$ under the null hypothesis. Then we have, as
$N\to\infty$,
\begin{equation}
\label{asymp_dist_MLE}
N^{1/3}\bigl\{\hat F_N(t_0)-F(t_0)\bigr\}\big/
\bigl\{\tfrac12F(t_0)(1-F(t_0))f(t_0)/g(t_0)\bigr\}^{1/3}
\stackrel{{\cal D}}{\longrightarrow}2Z,
\end{equation}
where $\stackrel{{\cal D}}{\longrightarrow}$ denotes convergence in distribution, and where
$Z$ is the last time where standard two-sided Brownian motion plus the parabola
$y(t)=t^2$ reaches its minimum.
\end{theorem}

From this one can deduce, under the assumptions of Theorem \ref{th:local_asymp_curstat_MSLE},
\begin{align}
\label{asymptotic_risk}
& N^{1/3} E \int_a^b\frac{N^{2/3}\bigl\{\hat F_N(x)-F(x)\bigr\}^2}{F(x)\{1-F(x)\}}\,dG(x)\sim N^{1/3}4EZ^2\int_a^b \frac{\bigl\{f(x)^2g(x)\bigr\}^{1/3}}{\bigl(4F(x)\{1-F(x)\}\bigr)^{1/3}}\,dx,\,N\to\infty,
\end{align}
where $Z$ is defined as in Theorem \ref{th:cur_stat_convergence}. By Table 4 in \shortciteN{piet_jon:01} we have:
$$
4EZ^2\approx1.05423856.
$$

Let $K_N$ be the number of jumps of the MLE on the interval $[a,b]$. Then it follows from \shortciteN{piet:11} that, again under the assumptions of Theorem \ref{th:local_asymp_curstat_MSLE},
\begin{equation}
\label{intensity_constant}
E K_n\sim cN^{1/3}\int_a^b \frac{\bigl\{f(x)^2g(x)\bigr\}^{1/3}}{\bigl(4F(x)\{1-F(x)\}\bigr)^{1/3}}\,dx,\,n\to\infty.
\end{equation}
for a constant $c>0$ which is close to $2.1$, so we find
$$
\frac{4EZ^2}{c}\approx 0.5
$$
It is tempting to believe that this ratio is exactly equal to $1/2$, but we have no proof of that.
It can also be deduced from \shortciteN{piet:11} that $K_N$ is asymptotically normal and that, in fact,
\begin{align}
\label{asymp_normality_jumps}
\frac{K_N-EK_N}{\sqrt{EK_N}}\stackrel{{\cal D}}\longrightarrow N(0,c_2),
\end{align}
for a universal constant $c_2>0$, not depending on the underlying distributions.

The intuitive interpretation of all this is that we have histograms with a random number of cells, where, under the null hypothesis ${\cal H}_0$, the number of cells has an asymptotic expectation which is proportional to the asymptotic expectation on the right-hand side of (\ref{asymptotic_risk}). Note that
\begin{align*}
&\sqrt{K_N}\left\{\frac{2T_N}{K_N}-\frac{4EZ^2}{c}\right\}=\sqrt{E K_N}\left\{\frac{2T_N}{K_N}-\frac{4EZ^2}{c}\right\}+o_p(1),
\end{align*}
and that
\begin{align*}
\sqrt{E K_N}\left\{\frac{2T_N}{K_N}-\frac{4EZ^2}{c}\right\}
=\sqrt{E K_N}\left\{\frac{2T_N}{EK_N}-1\right\}+\frac{4EZ^2}{c}\frac{EK_N-K_N}{\sqrt{EK_N}}+o_p(1),
\end{align*}
where $c$ is as in (\ref{intensity_constant}). Since
$$
\frac{EK_N-K_N}{\sqrt{EK_N}}\stackrel{{\cal D}}\longrightarrow N(0,c_2),
$$
where $c_2$ is defined as in (\ref{asymp_normality_jumps}), it is clear that $\sqrt{K_N}\left\{2T_N/K_N-1\right\}$ is an asymptotic pivot under ${\cal H}_0$ if and only if $\sqrt{EK_N}\left\{2T_N/EK_N-1\right\}$ is an asymptotic pivot under ${\cal H}_0$.

So the situation is somewhat similar to the situation in section \ref{section:MSLE}, but on the other hand much more complicated because of the fact that the MLEs are in fact histogram-type estimators, where the number of cells of the histograms is random, and because of the fact that the estimators $\hat F_{N1}$, $\hat F_{N2}$,  and $\hat F_N$ are nonlinear estimators which are also asymptotically nonlinear, which leads to non-standard limit distributions of the pointwise estimators $\hat F_{Ni}(t)$ and $\hat F_N(t)$, in contrast with the MSLEs $\tilde F_{Ni}(t)$ and $\tilde F_N(t)$ which have normal limit distributions. Another complication is that $\hat F_N$, $\hat F_{N1}$ and $\hat F_{N2}$ have jumps at different locations.

Nevertheless we want to include this original LR test in our comparisons and we use the bootstrap method of section \ref{section:crit_value_bootstrap} for generating critical values for this test.

\section{A bootstrap method for determining the critical value}
\label{section:crit_value_bootstrap}
\setcounter{equation}{0}

We propose the following method for determining the critical value for testing the null hypothesis that the two samples come from the same distribution for the likelihood ratio test, discussed in section \ref{section:MSLE}.

First compute a MSLE $\tilde F_{N,\tilde b_N}$ for the combined sample as discussed in section \ref{section:MSLE}, using a bandwidth $\tilde b_N\asymp N^{-1/5}$. Then, using the observations $T_1,\dots,T_m$ and $T_{m+1},\dots,T_N$ of the two samples, generate corresponding bootstrap values $\dd_1^*,\dots,\dd_m^*$ and $\dd_{m+1}^*,\dots,\dd_N^*$ by letting the $\dd_i^*$ be independent Bernoulli $(\tilde F_{N,\tilde b_N}(T_i))$ random variables. So in practice we generate quasi-random independent Uniform$(0,1)$ variables $U_i^*$ by using a random number generator, and let $\dd_i^*$ be equal to $1$ if $U_i^*<\tilde F_{N,\tilde b_N}(T_i)$ and zero otherwise. If the observation distributions, generating $T_1,\dots,T_m$ and $T_{m+1},\dots,T_N$, respectively, are different, this structure is preserved in this procedure; in the computation of the MSLEs $\tilde F_{Nj}^*$ in the bootstrap samples the estimates $\tilde g_{Nj}$ of $g_j$ in the original samples are used, for $j=1,2$.
Repeating this procedure $B$ times, we obtain $B$ bootstrap values $V^*_{N,i}$, $1\le i\le B$, of the test statistic. The distribution of $V_N$ under the null hypothesis is now approximated by the empirical distribution of these bootstrap values and the critical value at (for example) level $5\%$ by the $95$th percentile of this set of bootstrap values $V^*_{N,i}$.

In justifying this method for our test statistic $V_N$, we use the following theorem.

\begin{theorem}
\label{th:bootstrap1}
Let, under either of the conditions of Theorems \ref{th:local_asymp_curstat_MSLE} to \ref{th:SLR_test_same_g}, $\tilde{F}_{N,\tilde b_N}$ be the MSLE of $F$ under the null hypothesis, defined by the slope of the cusum diagram (\ref{cusum_combined_smooth}), where the bandwidth $\tilde b_N$ satisfies $\tilde b_N\asymp N^{-1/5}$.  Let $V_N^*$ be defined by
\begin{align}
\label{hist_LR_statistic*}
&V_N^*=\frac{2m}{N}\int_{t\in[a,b]} \left\{\tilde h_{N1}^*(t)\log \frac{{\tilde F}_{N1}^*(t)}{{\tilde F}_N^*(t)}+\bigl\{\tilde g_{N1}(t)-\tilde h_{N1}^*(t)\bigr\}\log\frac{1- {\tilde F}_{N1}^*(t)}{1- {\tilde F}_N^*(t)}\right\}\,dt\nonumber\\
&\qquad\qquad+\frac{2n}{N}\int_{t\in[a,b]} \left\{\tilde h_{N2}^*(t)\log \frac{{\tilde F}_{N2}^*(t)}{\tilde F_N^*(t)}+\bigl\{\tilde g_{N2}(t)-\tilde h_{N2}^*(t)\bigr\}\log\frac{1- {\tilde F}_{N2}^*(T_i)}{1- {\tilde F}_N^*(T_i)}\right\}\,dt
\end{align}
where $\tilde F_N^*$, $\tilde F_{N1}^*$ and $\tilde F_{N2}^*$ are the MSLEs, computed for the samples $(T_1,\dd_1^*),\dots,(T_m,\dd_m^*)$ and $(T_{m+1},\dd_{m+1}^*),\dots,(T_N,\dd_N^*)$,
and where the $\dd_i^*$ are Bernoulli $(\tilde F_{N,b_N}(T_i))$  random variables, generated in the way described before the statement of this theorem; $\tilde g_{Ni}$ and $\tilde h_{Ni}^*$ are  kernel estimates of the relevant observation densities, just as in section \ref{section:MSLE}, where
$$
\tilde h_{N1}^*(t)=m^{-1}\sum_{i=1}^m \dd_i^*K_{b_N}(t-T_i),\qquad\tilde h_{N2}^*(t)=n^{-1}\sum_{i=m+1}^N \dd_i^*K_{b_N}(t-T_i).
$$
with the same bandwidth $b_N$ as taken in the original samples, and where the densities $\tilde g_{N1}$ and $\tilde g_{N2}$ are the same as in the original samples.

Then we get under ${\cal H}_0$ that the conditional distribution function of $V_N^*$, given $(T_1,\dd_1)$, $\dots$, $(T_N,\dd_N)$, rescaled in the same way as in Theorems \ref{th:local_asymp_curstat_MSLE} to \ref{th:SLR_test_same_g} (depending on the choice of bandwidth $b_N$ and presence or absence of the condition $g_1=g_2$), converges at each $x\in\R$ in probability to the standard normal distribution function $\Phi(x)$.
\end{theorem}

The proof of this result is given in the appendix. If the null hypothesis does not hold, we follow the same scheme. The critical value is again determined by first computing the $\dd_i^*$, using the MSLE $\tilde F_{N,\tilde b_N}$, based on the combined sample.

For the MLEs of section \ref{section:LR} we follow a similar procedure, although we presently cannot justify this with a result analogous to Theorem \ref{th:bootstrap1}. However, the $\dd_i^*$'s are computed by using the MSLE $\tilde F_{N,\tilde b_N}$, based on the original combined sample, using a bandwidth $\tilde b_N\asymp N^{-1/5}$, instead of the ordinary MLE for this sample. This seems to work better for the sample sizes we used in the simulations. For these distributions, the MSLE converges at the local rate $N^{-2/5}$ instead of MLE itself, which has local rate $N^{-1/3}$, and this led to a better estimate of the level under the null hypothesis, which was taken to be $0.05$. Bootstrap estimates, based on the MLE instead of the MSLE, which we also computed, exhibited a very anti-conservative behavior for certain combinations of the parameters, sometimes leading to estimates of the levels which were twice the intended level.

\section{Other nonparametric tests}
\label{section:other_proposals}
\setcounter{equation}{0}
Most test which have been proposed for this problem are based on a comparison of simple functionals of the $\dd_i$. Under the assumption that the observation times $T_i$ have the same distribution in the two samples, the following test statistic is proposed in \shortciteN{sun:06}:
\begin{equation}
\label{suns_proposal1}
U_{cw}=\sum_{i=1}^N\left(Z_i-\bar Z\right)\dd_i=\b_N\sum_{i=1}^m \dd_i-\a_N\sum_{i=m+1}^N\dd_i,\qquad \a_N=\frac{m}{N},\qquad\b_N=\frac{n}{N},
\end{equation}
where we take $Z_i=1$ if the observation belongs to the first sample and $Z_i=0$ if the observation belongs to the second sample in the notation of \shortciteN{sun:06}, p.\ 76, and where $\bar Z=\sum_{i=1}^N Z_i/N$. 

It is stated in \shortciteN{sun:06} that the variance of $N^{-1/2}$ times (\ref{suns_proposal1}) is given by the random variable
\begin{equation}
\label{Suns_variance}
N^{-1}\left\{\sum_{i=1}^m \b_N^2\dd_i^2+\sum_{i=m+1}^N \a_N^2\dd_i^2\right\}.
\end{equation}
Apart from the fact that the variance then is a random variable, we have more difficulties in interpreting this, since we get, if $\a_N\to\a\in(0,1)$ and $\b_N\to\b=1-\a$,
\begin{align*}
&N^{-1}\left\{\sum_{i=1}^m \b_N^2\dd_i^2+\sum_{i=m+1}^N \a_N^2\dd_i^2\right\}
\stackrel{p}\longrightarrow \a\b\left\{\b\int F(t)\,dG_1(t)+\a\int F(t)\,dG_2(t)\right\}\\
&=\a\b\int F(t)\,dG(t),
\end{align*}
if $G_1=G_2=G$. But the actual variance of $N^{-1/2}$ times (\ref{suns_proposal1}) is given by:
\begin{equation}
\label{actual_variance}
\a_N\b_N\int F(t)\,dG(t)\left\{1-
\int F(t)\,dG(t)\right\},
\end{equation}
if $G_1=G_2=G$. So the proposed estimate of the variance in \shortciteN{sun:06} will severely overestimate the actual variance, and the proposed normalization will not give a standard normal distribution in the limit, as claimed in \shortciteN{sun:06}.

Also, considering the $Z_i$ as i.i.d.\ random variables, as in \shortciteN{sun:06}, where $Z_i$ is a Bernoulli random variable with
$$
\P\{Z_i=1\}=\a_N,\qquad \P\{Z_i=0\}=\b_N=1-\a_N,
$$
and where the $Z_i$ are independent of the observation times $T_i$ and the indicators $\dd_i$, we arrive at (\ref{actual_variance}) instead of (\ref{Suns_variance}) as the approximate variance of $N^{-1/2}$ times (\ref{suns_proposal1}). This is seen in the following way.

We can write, under the null hypothesis that the $\dd_i$ have the same distribution, and also under the restriction that the observations $T_i$ have the same distribution in the two samples,
\begin{align}
\label{decomposition_Sun_statistic}
U_{cw}&=\sum_{i=1}^N\left(Z_i-\bar Z\right)\dd_i=\sum_{i=1}^N\left(Z_i-\bar Z\right)\left(\dd_i-E\dd_i\right)\nonumber\\
&=\sum_{i=1}^N\left(Z_i-\a_N\right)\left(\dd_i-E\dd_1\right)+\left(\a_n-\bar Z\right)\sum_{i=1}^N\left(\dd_i-E\dd_1\right),
\end{align}
using $E\dd_i=E\dd_1$ for each $i$. This yields:
\begin{align*}
N^{-1}\mbox{var}\left(U_{cw}\right)\sim \mbox{var}\left(\left(Z_1-\a_N\right)\left(\dd_1-E\dd_1\right)\right)
=\a_N\b_N\int F(t)\,dG(t)\left\{1-
\int F(t)\,dG(t)\right\},
\end{align*}
since the second expression on the right-hand side of (\ref{decomposition_Sun_statistic}) gives a contribution of lower order. So we arrive (not surprisingly) again at (\ref{actual_variance}) as an approximation of the variance of $n^{-1/2}U_{cw}$ in the interpretation of the $Z_i$ as i.i.d.\ random variables, implying that the $\hat \s_{cw}$ suggested as standardization of the statistic $N^{-1/2}U_{cw}$ in \shortciteN{sun:06} in the last line of the first paragraph of section 4.2.1.1, has to be replaced by an estimate of the square root of (\ref{actual_variance}), also if we consider the $Z_i$ to be random. The mistake of taking (\ref{Suns_variance}) as an estimate of the variance is probably caused by ignoring the dependence of the terms $(Z_i-\bar Z)\dd_i$, caused by $\bar Z$, and treating $\bar Z$ as if it were $EZ_i$. The presence of $\bar Z$ actually has a variance diminishing effect.

Putting these difficulties aside, and not using the standardization by the square root of (\ref{Suns_variance}), we could of course consider the test statistic
\begin{equation}
\label{test_Sun1}
\tilde U_N=N^{-1/2}\left\{\b_N\sum_{i=1}^m \dd_i-\a_N\sum_{i=m+1}^N\dd_i,\right\}
\end{equation}
which has expectation zero under the null hypothesis, provided $G_1=G_2$, and variance (\ref{actual_variance}), if $G_1=G_2=G$. Then, since the MLE $\hat F_N$, based on the combined samples, satisfies, under some regularity conditions,
\begin{align*}
&\int \hat F_N(t)\,d\G_N(t)\stackrel{p}\longrightarrow \int F(t)\,dG(t),
\end{align*}
where $F$ is the limit (mixture) distribution of the combined samples (which is the underlying distribution under ${\cal H}_0$), we could use as test statistic
\begin{equation}
\label{test_Sun}
U_N=\frac{\tilde U_N}{\hat\s_N}.
\end{equation}
where $\tilde U_n$ is defined by (\ref{test_Sun1}), and where
\begin{equation}
\label{variance_Sun_corrected}
\hat\s_N^2=\a_N\b_N\int \hat F_N(t)\,d\G_N(t)\left\{1-\int \hat F_N(t)\,d\G_N(t)\right\},
\end{equation}
Then $U_N$ tends to a standard normal distribution under the null hypothesis, if $G_1=G_2=G$.
We note that in \shortciteN{sun:06} also a test where $G_1\ne G_2$ is allowed is discussed, but since this test is connected to a specific parametric model, it is not a test of the fully nonparametric type we consider here.

\shortciteN{andersen_ronn:95} consider a test based on
$$
W_N=\frac{\sqrt{N}\int_0^a\bigl\{\hat F_{N1}(t)^2-\hat F_{N2}(t)^2\bigr\}\,d\G_N(t)}
{\sqrt{\frac4{\a_N\b_N}\int_0^a \hat F_N(t)^3\bigl\{1-\hat F_N(t)\bigr\}\,d\G_N(t)}}\,,
$$
on an interval $[0,a]$, where $W_N$ is asymptotically standard normal under the null hypothesis, if $G_1=G_2$ (note that in their definition of this test statistic, which is denoted by $W$ on p.\ 325, a factor $\sqrt{n}$ is missing in the numerator). They rely in their proof on the master's thesis \shortciteN{bettina:91}, which, incidentally, was written at Delft University of Technology, and not at the University of Copenhagen, as stated in \shortciteN{andersen_ronn:95}.

Under the conditions of Theorem \ref{th:local_asymp_curstat_MSLE} we have:
\begin{align}
\label{Andersen_asymp}
\frac{\sqrt{N}\int_{[a,b]}\bigl\{\hat F_{N1}(t)^2-\hat F_{N2}(t)^2\bigr\}\,d\G_N(t)}
{\sqrt{\frac4{\a_N\b_N}\int_{[a,b]} \hat F_N(t)^3\bigl\{1-\hat F_N(t)\bigr\}\,d\G_N(t)}}
\stackrel{{\cal D}}\longrightarrow N(0,1),
\end{align}
under ${\cal H}_0$, where $N(0,1)$ is the standard normal distribution.
A sketch of how this result can be derived, roughly using the techniques developed in \shortciteN{bettina:91}, is given in the appendix. 

\section{A simulation study}
\label{section:sim_study}
\setcounter{equation}{0}
In this section we compare the LR test based on the MSLEs and the real LR test with the methods, discussed in the preceding section. In our comparison we use the same Weibull model, which was used in the comparison, given in \shortciteN{andersen_ronn:95}. In determining the critical levels and the powers of the tests, based on $V_N$ (the test statistic based on the MSLEs) and the LR test, based on the MLEs, we used the method described in section \ref{section:crit_value_bootstrap}, that is, the critical values were determined by (Bernoulli) bootstrapping the $\dd_i$, using the MSLE $\tilde F_{N,\tilde b_N}(T_i)$ for the combined samples at the observations $T_i$, by taking $1000$ bootstrap samples and determining the $95$th percentile of the bootstrap test statistics, so obtained.

As the bandwidth for smoothing the MLE $\hat F_N$, we used $b_N=2N^{-1/5}$ in all instances, and we used the kernel (\ref{triweight}) in computing $\tilde F_N$, as described in section \ref{section:MSLE}.
As the observation densities $g_1$ and $g_2$ for the observation times $T_i$ we took the uniform densities on $[0,2]$, just as in \shortciteN{andersen_ronn:95}. Note that in the simulation study of \shortciteN{andersen_ronn:95} $g_1=g_1$, so we can apply Theorem \ref{th:SLR_test_same_g}. This allowed us to resample from the MSLE $\tilde F_N$, which was also used in the computation of the test statistic for the original samples. 

The powers and levels computed below for the test statistics $V_N$ (MSLEs) and the LR statistic, based on the MLEs, are determined by taking $1000$ samples from the original distributions and taking $1000$ bootstrap sample from each sample, rejecting the null hypothesis if the value in the original sample was larger than the $950$th order statistic of the values obtained in the bootstrap samples. The values given in the tables below represent the fraction of rejections for the $1000$ samples from the original distributions. The simulation were carried out using a $C$ program, which was written by the author specifically for this analysis.

We also included the estimates, discussed in section \ref{section:other_proposals}, where $W_N$ denotes the test statistic of \shortciteN{andersen_ronn:95} and $U_N$ denotes the test statistic of \shortciteN{sun:06}, but with the incorrect estimate of the variance (\ref{Suns_variance}) in \shortciteN{sun:06} replaced by (\ref{variance_Sun_corrected}). In this case we just took $1.96$ as our critical value for the absolute value of the test statistic, since the convergence to the standard normal distribution is reasonably fast for these test statistics under the null hypothesis. In this way one can rather fastly compute tables of this type for these test statistics, which was again done by writing a $C$ program for this purpose. The tabled values are again based on $1000$ samples from the original (Weibull) distributions.

Using the same parametrization as in \shortciteN{andersen_ronn:95}, we generated the first sample from the density
\begin{equation}
\label{Weibull1}
\a_1\l x^{\a_1-1} e^{-\l x^{\a_1}},\,x>0,
\end{equation}
and the second sample from the density
\begin{equation}
\label{Weibull2}
\a_2\l\theta x^{\a_2-1} e^{-\l\theta x^{\a_2}},\,x>0,
\end{equation}
where $\l=1.6$ or $\l=0.58$, and $\a_1=0.5, 1.0$ or $2.0$. The value of $\th$ is $1, 1.25$ or $2$. Why these specific values were taken in \shortciteN{andersen_ronn:95} is not clear to me, but I take the same values for an easy comparison with the work, reported in their paper. I have to note, though, that for $\a_i=0.5$ the Weibull density is unbounded near zero, and that then the results of \shortciteN{bettina:91} are not valid on $[0,2]$, since one of the conditions in her thesis  was that this density is bounded on the interval of interest. This is also one of the reasons that the interval $[0,2]$, used in \shortciteN{andersen_ronn:95}, was shrunk to $[0.1,1.9]$ in our simulation study, since the density is bounded on this interval.

To illustrate the effect of different observation distributions in the two samples, we generated the first sample of $T_i$'s again from the uniform density on $[0,2]$, but the second sample from the decreasing density
$$
g_2(t)=\tfrac14(2-t)^3,\,t\in[0,2],
$$
see Tables \ref{table:levels1_different_g} and \ref{table:levels2_different_g}. Note that in this case Theorem \ref{th:SLR_test_same_g} does not apply, and we would actually have to use Theorem \ref{th:local_asymp_curstat_MSLE} or \ref{th:SLR_test}. Nevertheless, we just proceeded in the same way as for the simulations for the situation $g_1=g_2$, and Tables \ref{table:levels1_different_g} and \ref{table:levels2_different_g} show that the test based on the MSLEs, where we take $b_N=2N^{-1/5}$ and compute the critical values using the bootstrap procedure, were rather insensitive to the difference of the observation distributions $G_1$ and $G_2$.

\begin{table}[ht!]
\centering
\caption{Estimated levels. The estimation interval is $[0.1,1.9]$, and $m=n=50$; $g_1(t)\equiv\tfrac12$, $g_2(t)\equiv\tfrac12$, $\a_1=\a_2$. The intended level is $\a=0.05$.}
\label{table:levels1}
\vspace{0.5cm}
\begin{tabular}{|l|c|c|c|c|c|c|c|c|c|c|c|c|}
\hline
 $g_1= g_2$ & $\lambda,\a_i$ & \multicolumn{5}{c|}{Under $H_0$}\\
\hline
$m=n=50$  & $1.6,0.5$ &  $1.6,1.0$ & $1.6,2.0$ & $0.58, 0.5$ & $0.58, 1.0$ & $0.58, 2.0$\\
\hline
SLR test & 0.041
 & 0.058 & 0.045 & 0.049 & 0.049 &0.059\\
LR test  &0.045 & 0.051 & 0.041  & 0.052 & 0.046 & 0.055\\
$U_N$ & 0.050 & 0.060& 0.047& 0.054 & 0.058 & 0.052  \\
$W_N$ &0.055  & 0.066 & 0.087 & 0.061 & 0.061 & 0.072\\
\hline
\end{tabular}
\end{table}

\begin{table}[ht!]
\centering
\caption{Estimated levels. The estimation interval is $[0.1,1.9]$, and $m=n=50$; $g_1(t)\equiv\tfrac12$, $g_2(t)=\tfrac14(2-t)^3$, $\a_1=\a_2$. The intended level is $\a=0.05$.}
\label{table:levels1_different_g}
\vspace{0.5cm}
\begin{tabular}{|l|c|c|c|c|c|c|c|c|c|c|c|c|}
\hline
$g_2(t)=\tfrac14(2-t)^3$ & $\lambda,\a_i$ & \multicolumn{5}{c|}{Under $H_0$}\\
\hline
$m=n=50$  & $1.6,0.5$ &  $1.6,1.0$ & $1.6,2.0$ & $0.58, 0.5$ & $0.58, 1.0$ & $0.58, 2.0$\\
\hline
SLR test & 0.049
 & 0.051 & 0.045 & 0.049 & 0.049 &0.059\\
LR test  &0.051 & 0.055 & 0.049  & 0.044 & 0.050 & 0.056\\
$U_N$ & 0.422 & 0.745& 0.950& 0.262 & 0.540 & 0.885  \\
$W_N$ &0.122  & 0.108 & 0.130 & 0.326 & 0.302 & 0.276\\
\hline
\end{tabular}
\end{table}

\begin{table}[ht!]
\centering
\caption{Estimated levels. The estimation interval is $[0.1,1.9]$, and $m=n=250$; $g_1(t)\equiv\tfrac12$, $g_2(t)\equiv\tfrac12$, $\a_1=\a_2$. The intended level is $\a=0.05$.}
\label{table:levels2}
\vspace{0.5cm}
\begin{tabular}{|l|c|c|c|c|c|c|c|c|c|c|c|c|}
\hline
 $g_1= g_2$, & $\lambda,\a_i$ & \multicolumn{5}{c|}{Under $H_0$}\\
\hline
$m=n=250$  & $1.6,0.5$ &  $1.6,1.0$ & $1.6,2.0$ & $0.58, 0.5$ & $0.58, 1.0$ & $0.58, 2.0$\\
\hline
SLR test & 0.051 & 0.049 & 0.052 & 0.053 & 0.032 &0.040\\
LR test  &0.048 &  0.049 & 0.059  & 0.053 & 0.045 & 0.054\\
$U_N$ & 0.050 & 0.060& 0.047& 0.054 & 0.058 & 0.052  \\
$W_N$ &0.055 & 0.066 & 0.087 & 0.061 & 0.061 &0.072  \\
\hline
\end{tabular}
\end{table}

\begin{table}[ht!]
\centering
\caption{Estimated levels. The estimation interval is $[0.1,1.9]$, and $m=n=250$. The intended level is $\a=0.05$; $g_1(t)=\tfrac12$, $g_2(t)=\tfrac14(2-t)^3$, $\a_1=\a_2$.}
\label{table:levels2_different_g}
\vspace{0.5cm}
\begin{tabular}{|l|c|c|c|c|c|c|c|c|c|c|c|c|}
\hline
 $g_2(t)=\tfrac14(2-t)^3$& $\lambda,\a_i$ & \multicolumn{5}{c|}{Under $H_0$}\\
\hline
$m=n=250$  & $1.6,0.5$ &  $1.6,1.0$ & $1.6,2.0$ & $0.58, 0.5$ & $0.58, 1.0$ & $0.58, 2.0$\\
\hline
SLR test & 0.044 & 0.050 & 0.051 & 0.049 & 0.044 &0.051\\
LR test  &0.045 & 0.051 & 0.041  & 0.052 & 0.054 & 0.058\\
$U_N$ & 0.970 & 1.000& 1.000& 0.840 & 0.996 & 1.000  \\
$W_N$ & 0.181 & 0.135& 0.102& 0.513 & 0.491 & 0.410\\
\hline
\end{tabular}
\end{table}

\begin{table}[ht!]
\centering
\caption{Powers for different shapes, if $m=n=50$. The estimation interval is $[0.1,1.9]$.}
\label{table:powers2}
\vspace{0.5cm}
\begin{tabular}{|l|c|c|c|c|c|c|c|c|c|c|c|c|}
\hline
 $g_1=g_2$ & $\lambda,\a_1,\a_2$ & \multicolumn{3}{c|}{Different shapes}\\
\hline
$m=n=50$ & $1.6,0.5, 1.0$ &  $1.6,0.5, 2.0$ & $0.58, 0.5 ,2.0$ & $0.58, 1.0, 2.0$ \\
\hline
SLR test & 0.174 & 0.675 & 0.470 & 0.207 \\
LR test &0.125 & 0.533 & 0.364  & 0.173 \\
$U_N$ & 0.061& 0.069 & 0.045& 0.053   \\
$W_N$ &0.062 & 0.110 & 0.179 & 0.146 \\
\hline
\end{tabular}
\end{table}

\begin{table}[ht!]
\centering
\caption{Powers for different shapes, if $m=n=250$. The estimation interval is $[0.1,1.9]$.}
\label{table:powers3}
\vspace{0.5cm}
\begin{tabular}{|l|c|c|c|c|c|c|c|c|c|c|c|c|}
\hline
 $g_1=g_2$ & $\lambda,\a_1,\a_2$ & \multicolumn{3}{c|} {Different shapes}\\
\hline
$m=n=250$ & $1.6,0.5, 1.0$ &  $1.6,0.5, 2.0$ & $0.58, 0.5 ,2.0$ & $0.58, 1.0, 2.0$ \\
\hline
SLR test & 0.606 & 1.000 & 0.990 & 0.787 \\
LR test &0.440 & 1.000 & 0.974  & 0.610 \\
$U_N$ & 0.076& 0.132 & 0.062& 0.076  \\
$W_N$ &0.088 & 0.112 & 0.583 & 0.406 \\
\hline
\end{tabular}
\end{table}

\begin{table}[ht!]
\centering
\caption{Powers for different baseline hazards, same shape, if $m=n=50$. The estimation interval is $[0.1,1.9]$. The parameters $\a_i$ are either both $0.5$ or both $2$ and $\l=1.6$ or $0.58$; $\theta=1.25, 1.5$ or $2$.}
\label{table:powers4}
\vspace{0.5cm}
\begin{tabular}{|l|c|c|c|c|c|c|c|c|c|c|c|c|}
\hline
 $g_1=g_2$ & $\lambda,\a_i,\th$ & \multicolumn{5}{c|} {Different baseline hazards}\\
\hline
$m=n=50$ & $1.6,0.5, 1.25$ &  $1.6, 0.5, 1.5$ & $1.6, 0.5, 2$ & $0.58, 2, 1.25$ & $0.58, 2, 1.5$ & $0.58, 2, 2$\\
\hline
SLR test  & 0.138 & 0.283 & 0.632 & 0.091 & 0.208 &0.480\\
LR test &0.097 & 0.218 & 0.498  & 0.082 & 0.171 & 0.342\\
$U_N$ & 0.108 & 0.198& 0.441& 0.100 & 0.151 & 0.333  \\
$W_N$ &0.147 & 0.352 & 1.000 & 0.103 & 0.293 & 0.681\\
\hline
\end{tabular}
\end{table}

\begin{table}[ht!]
\centering
\caption{Powers for different baseline hazards, same shape, if $m=n=250$. The estimation interval is $[0.1,1.9]$. The parameters $\a_i$ are either both $0.5$ or both $2$ and $\l=1.6$ or $0.58$; $\th=1.25, 1.5$ or $2$.}
\label{table:powers5}
\vspace{0.5cm}
\begin{tabular}{|l|c|c|c|c|c|c|c|c|c|c|c|c|}
\hline
 $g_1=g_2$ & $\lambda,\a_i,\th$ & \multicolumn{5}{c|} {Different baseline hazards}\\
\hline
$m=n=250$ & $1.6,0.5, 1.25$ &  $1.6, 0.5, 1.5$ & $1.6, 0.5, 2$ & $0.58, 2, 1.25$ & $0.58, 2, 1.5$ & $0.58, 2, 2$\\
\hline
SLR test  & 0.377 & 0.873 & 1.000 & 0.227 & 0.689 &0.995\\
LR test &0.246 & 0.728 & 0.996  & 0.171 & 0.505 & 0.964\\
$U_N$ & 0.324 & 0.721& 0.971& 0.200 & 0.495 & 0.921  \\
$W_N$ &0.473 & 0.912 & 1.000 & 0.337 & 0.835 & 1.000\\
\hline
\end{tabular}
\end{table}

The results of our experiments can be summarized in the following way. The corrected version of the test statistic discussed in \shortciteN{sun:06}, denoted by $U_N$ here, has almost no power for different shape alternatives of the type shown in Figure \ref{fig:MSLE+MLE_Weibull}, even for sample sizes $m=n=250$. The test proposed by \shortciteN{andersen_ronn:95}, denoted by $W_N$, has somewhat more power here, but is clearly also not very good for this type of alternative, as already discussed in \shortciteN{andersen_ronn:95} (they call this the ``crossing alternatives", since the distribution functions indeed cross). Both the test based on the MSLEs and the test, based on the MLEs, have more power here.
The test, based on $W_N$, is surprisingly powerful for the alternatives which have the same shape but different baseline hazards, and the test, based on $U_N$ also has more power here. The other tests, based on the MSLEs and MLEs, have also reasonable power here, in particular the test based on the MSLEs. Finally, Tables \ref{table:levels1_different_g} and \ref{table:levels2_different_g} show that the observation distributions in the two samples can be different if we use the LR-type tests, in contrast with the other tests, considered here. In fact, it has a disastrous effect for the tests $U_N$ and $W_N$; $U_N$ even gives $100\%$ rejection under the null hypothesis for several combinations of the parameters.

As noted in the introduction, one could try to use a permutation distribution approach in estimating the levels of the tests under the null hypothesis, also when the observation distributions are different. This does not seem to make much sense for the tests, based on $U_N$ and $W_N$, but could possibly be of use for the tests, based on the MSLEs and MLEs. We did some experiments in this direction for the Weibull distributions of the simulation study, with rather bad results for our sample sizes $m=n=50$ and $m=n=250$. The general finding is that the test based on the MLEs becomes very conservative, whereas the estimates of the levels for the tests based on the MSLEs become too variable to be of any use. In the latter case one big difference with the approach using the bootstrapped $\dd_i$ is that for the approach using the permutation distribution, the densities $g_1$ and $g_2$ have to be estimated anew for every new permutation of the variables $(T_1,\dd_1),\dots,(T_N,\dd_N)$, whereas these estimates can be held fixed in the bootstrap approach. This probably leads to a higher variability of the values of the test statistic under the null hypothesis for the permutation approach, leading to unstable estimates of the levels.
However, when the observation distributions are the same in the two samples, the permutation procedure seems to work fine, and then gives the same results as the bootstrap procedure.

As a general rule one can say that the tests, based on $U_N$ or $W_N$, can only have power if the corresponding moment functionals are different from zero. For $U_N$ this functional is given by
\begin{equation}
\label{sun_functional}
\int_a^b\{F_1(t)-F_2(t)\}\,dG(t),
\end{equation}
and for $W_N$ it is given by
\begin{equation}
\label{Andersen_functional}
\int_a^b\{F_1(t)^2-F_2(t)^2\}\,dG(t).
\end{equation}
It is clear that $F_1$ and $F_2$ can be very different and still satisfy
$$
\int_a^b\{F_1(t)-F_2(t)\}\,dG(t)=0,\qquad\mbox{ or }\qquad \int_a^b\{F_1(t)^2-F_2(t)^2\}\,dG(t)=0
$$
and in that case that tests, based on $U_N$ or $W_N$, respectively, will have no power. The LR tests will not suffer from this drawback, since they involve a Kullback-Leibler type distance, and are locally (for example if one would consider contiguous alternatives) equivalent to the squared $L_2$-distance
\begin{equation}
\label{ML_functional}
\int_a^b\frac{\{F_1(t)-F(t)\}^2}{F(t)\{1-F(t)\}}\,dG_1(t)+\int_a^b\frac{\{F_2(t)-F(t)\}^2}{F(t)\{1-F(t)\}}\,dG_2(t),
\end{equation}
where $F$ is the distribution function of the combined sample. Moreover, they allow the observation distributions to be different in the two samples, something the other test also do not allow.

\begin{figure}[!ht]
\begin{center}
\includegraphics{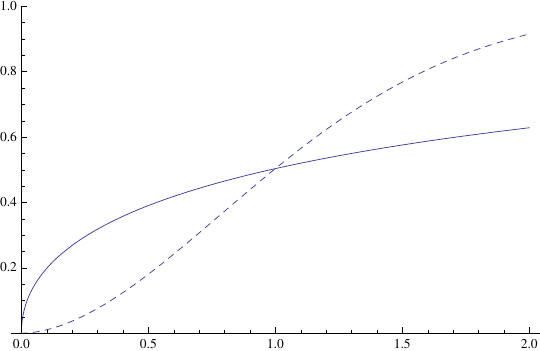}
\end{center}
\caption{The Weibull distribution function with parameters $\a=0.5$ and $\l=0.7$ (solid curve) and the Weibull distribution function with parameters $\a=1.8153$ and $\l=0.7$ (dashed).}
\label{fig:Weibull1}
\end{figure}

The Weibull alternatives, considered in the simulation study, form a family for which the integrals, corresponding to the statistics $U_N$ and $W_N$ are different under the alternatives, considered there. So for these type alternatives the tests $U_N$ and $W_N$ can be expected to have a power exceeding the level of the test. But if the first sample is generated from a Weibull distribution function $F_1$ with parameters $\a=0.5$ and $\l=0.7$ and the second sample is generated from a Weibull distribution function $F_2$ with parameters $\a=1.8153$ and $\l=0.7$, the distribution functions are very different (see Figure \ref{fig:Weibull1}), although we get:
$$
\int_a^b\left\{F_1(t)-F_2(t)\right\}\,dt\approx-1.87\cdot10^{-6},\qquad a=0.1,\qquad\,b=1.9.
$$
Taking again the observations $G_1$ and $G_2$ to be uniform on $[0,2]$, we get that the test based on the MSLE has power $0.993$ for this alternative, whereas the tests based on $U_N$ has power $0.048$ (which is lower than the level $0.05$).

\begin{figure}[!ht]
\begin{center}
\includegraphics{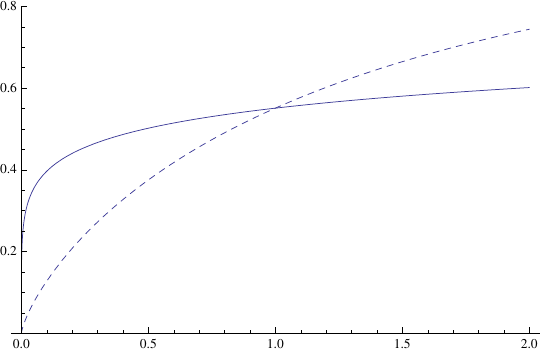}
\end{center}
\caption{The Weibull distribution function with parameters $\a=0.2$ and $\l=0.8$ (solid curve) and the Weibull distribution function with parameters $\a=0.767$ and $\l=0.8$ (dashed).}
\label{fig:Weibull2}
\end{figure}

If the first sample is generated from a Weibull distribution function $F_1$ with parameters $\a=0.2$ and $\l=0.8$ and the second sample is generated from a Weibull distribution function $F_2$ with parameters $\a=0.767$ and $\l=0.8$, the distribution functions are again rather different (see Figure \ref{fig:Weibull2}), although we get:
$$
\int_a^b\left\{F_1(t)^2-F_2(t)^2\right\}\,dt\approx2.6\cdot10^{-6},\qquad a=0.1,\qquad\,b=1.9.
$$
Taking $g_1=g_2\equiv(1/2)1_{[0,2]}$ again, the test based on the MSLE has power $0.713$ for this alternative, whereas the tests based on $W_N$ has power $0.041$ (which is again lower than the level $0.05$).

The LR tests, based on the MLEs instead of the MSLEs, has powers $0.964$ and $0.515$, respectively, for these alternatives, taking the sample sizes $m=n=250$ again.

\section{Concluding remarks}
\label{section:conclusion}
In the preceding, two fully nonparametric tests for the two-sample problem for current status data were discussed. The tests allow the observation distributions for the two samples to be different, and will be consistent for any situation where (\ref{ML_functional}) will be different from zero and the distributions satisfy some regularity conditions. For the test, based on the maximum smoothed likelihood estimators (MSLEs), the theory is more complete than for the test, based on the  MLEs, but we suggest a bootstrap method for determining critical values for the latter test, which seemed to work well in the simulation study we conducted.

Most tests which have been proposed for this problem rely on specific functionals, such as (\ref{sun_functional}) or (\ref{Andersen_functional}), which can easily be zero, while the distributions $F_1$ and $F_2$ are very different. If these functionals are zero, the tests cannot be expected to have power against these alternatives. A simulation study in section \ref{section:sim_study}, using a Weibull model, which was also used in \shortciteN{andersen_ronn:95}, further illustrates this point.

The convergence to normality in Theorems \ref{th:local_asymp_curstat_MSLE} to \ref{th:SLR_test_same_g} cannot be expected to be very fast. This phenomenon is well-known from the theory of integrated mean squared errors of density estimators. However, the bootstrap procedure we propose for estimating the critical values of the tests, discussed in section \ref{section:crit_value_bootstrap} seems to work well, even for sample sizes $m=n=50$. So, for practical purposes, we advise to use this procedure for estimating the critical values of the tests, instead of relying on the asymptotic normality under the null hypothesis.

We have chosen to work with conditional tests, and in this approach we only have to resample the $\dd_i$ in estimating the critical value for the tests. It is also possible to work with unconditional tests, but in that case one also has to resample the $T_i$ from estimates of the densities $g_1$ and $g_2$ for the first and second sample, respectively. Preliminary experiments with this procedure indicate that the resulting powers are roughly the same for the model, used in the simulation section \ref{section:sim_study}, but more research is needed to evaluate the two approaches. 

\newpage
\section{Appendix}
\label{section:appendix}
\setcounter{equation}{0}

\begin{lemma}
\label{Talagrand_lemma}
Let either of the conditions of Theorems \ref{th:local_asymp_curstat_MSLE} to \ref{th:SLR_test_same_g} be satisfied. Then
\begin{equation}
\label{sup_distance1}
\sup_{t\in[a,b]}\left|\tilde g_{Nj}(t)-g_j(t)\right|=O_p\left(N^{-(1-\a)/2}\sqrt{\log N}\right)
\end{equation}
and
\begin{equation}
\label{sup_distance2}
\sup_{t\in[a,b]}\left|\tilde h_{Nj}(t)-F(t)g_j(t)\right|=O_p\left(N^{-(1-\a)/2}\sqrt{\log N}\right),\,j=1,2,
\end{equation}
\end{lemma}
implying that also:
\begin{align*}
\sup_{t\in[a,b]}\left|\tilde F_{Nj}(t)-F(t)\right|=O_p\left(n^{-(1-\a)/2}\sqrt{\log N}\right),\,j=1,2.
\end{align*} 

\noindent
{\bf Proof.} By Corollary 3.4 in \shortciteN{piet_geurt_birgit:10} we have, with probability tending to one,
\begin{equation}
\label{rep_MSLE1}
\tilde F_{N1}(t)=\frac{\tilde h_{N1}(t)}{\tilde g_{N1}(t)},\,t\in[a,b],
\end{equation}
that is, the MSLE is just equal to the ratio of two kernel estimators for $t\in[a,b]$, with probability tending to one. Similarly, with probability tending to one,
\begin{equation}
\label{rep_MSLE2}
\tilde F_{N2}(t)=\frac{\tilde h_{N2}(t)}{\tilde g_{N2}(t)},\,t\in[a,b],
\end{equation}
and
\begin{equation}
\label{rep_MSLE_joint}
\tilde F_N(t)=\frac{\a_N\tilde h_{N1}(t)+\b_N\tilde h_{N2}(t)}{\a_N\tilde g_{N1}(t)+\b_N\tilde g_{N2}(t)},\,t\in[a,b],
\qquad\a_n=m/N,\,\qquad\b_N=1-\a_N.
\end{equation}
Hence we assume in the following that $\tilde F_N$, $\tilde F_{N1}$ and $\tilde F_{N2}$ have the representations (\ref{rep_MSLE1}), (\ref{rep_MSLE2}) and (\ref{rep_MSLE_joint}), respectively.

We consider the set of functions
\begin{equation}
\label{VC-class}
{\cal F}=\left\{\f:\f\bigl(x,u\bigm|t,h\bigr)=K\left(\frac{t-u}{h}\right)1_{[0,u]}(x),\,t\in[a,b],\,u\ge0,\,h\in(0,c]\right\},
\end{equation}
where $0<c\le(1/2)\min[a,M-b]$, where $M$ is the smallest number such that $\min\{F_1(M),F_2(M)\}=1$. The kernels, considered in this paper (see section \ref{section:MSLE}) satisfy the condition (K1) of \shortciteN{evariste:02}, p.\ 911, implying that $\cal F$ is a bounded VC class of measurable functions.
Furthermore,
\begin{align*}
&\mbox{var}\left(\f\bigl(X_1,T_1\bigm|t,h\bigr)\right)=\mbox{var}\left(K\left(\frac{t-T_1}{h}\right)\dd_1\right)=\int F(u)\{1-F(u)\}K\left(\frac{t-u}{h}\right)^2g_1(u)\,du\\
&=h\int F(t-hw)\{1-F(t-hw)\}K(w)^2g_1(t-hw)\,dw\le cK(0)^2\sup_{t\in[a-h,b+h]}g_1(t).
\end{align*}
Letting $\s$ and $U$ be defined as in Corollary 2.2 of \shortciteN{evariste:02}, we get from (2.8) in this corollary the following inequality, based on \shortciteN{talagrand:94} and \shortciteN{talagrand:96},
\begin{align}
\label{Talagrand_ineq}
&\P\left\{\sup_{t\in[a,b]}\left|\sum_{i=1}^m\left\{K\left(\frac{t-T_i}{h}\right)\dd_i-EK\left(\frac{t-T_i}{h}\right)\dd_i\right\}\right|\ge C\s\sqrt{m\log\left(\frac{K(0)}{\s}\right)}\right\}\nonumber\\
&\le L\exp\left\{-\frac{C\log\left\{1+C/(4L)\right\}}{L}\log(K(0)/\s)\right\},
\end{align}
where $L$ and $C$ are positive constants depending on the VC characteristics of the class ${\cal F}$, and where $\s$, specialized to our situation, is given by
$$
\s=K(0)\left(c\sup_{t\in[a-c,b+c]}g_1(t)\right)^{1/2}.
$$

Since we take the bandwidth $b_N$ of order $b_N\asymp N^{-\a}$, we get from (\ref{Talagrand_ineq}), taking $c$ in (\ref{VC-class}) also of order $O(n^{-\a})$,
\begin{align*}
&\sup_{t\in[a,b]}\left|m^{-1}\sum_{i=1}^m K_{b_N}\left(t-T_i\right)\dd_i-EK_{b_N}\left(t-T_1\right)\dd_1\right|=O_p\left(\sqrt{\frac1{Nb_n}\log\frac{K(0)}{\sqrt{b_N}}}\right)\\
&=O_p\left(N^{-(1-\a)/2}\sqrt{\log N}\right).
\end{align*}
Since we get directly from Theorem 2.3 in \shortciteN{evariste:02} that
$$
\sup_{t\in[a,b]}\left|\tilde g_{N1}(t)-EK_{b_N}\left(t-T_1\right)\right|=O_p\left(N^{-(1-\a)/2}\sqrt{\log N}\right).
$$
It now follows from (\ref{rep_MSLE1}), which holds with probability tending to one, that also
\begin{align*}
&\sup_{t\in[a,b]}\left|\tilde F_{N1}(t)-\frac{EK_{b_N}\left(t-T_1\right)\dd_1}{EK_{b_N}\left(t-T_1\right)}\right|=O_p\left(N^{-(1-\a)/2}\sqrt{\log N}\right).
\end{align*}
By the conditions of Theorem \ref{th:local_asymp_curstat_MSLE} we also have:
\begin{align*}
EK_{b_N}\left(t-T_1\right)\dd_1=\int K_{b_N}\left(t-u\right)F(u)g_1(u)\,du=F(t)g_1(t)+O\left(N^{-2\a}\right),
\end{align*}
and
\begin{align*}
EK_{b_N}\left(t-T_1\right)=g_1(t)+O\left(N^{-2\a}\right),
\end{align*}
uniformly for $t\in[a,b]$. Hence we obtain:
\begin{align*}
\sup_{t\in[a,b]}\left|\tilde F_{N1}(t)-F(t)\right|=O_p\left(n^{-(1-\a)/2}\sqrt{\log N}\right).
\end{align*}
The other relations are proved in a similar way.\eop

\begin{lemma}
\label{expectation_lemma}
Let either of the conditions of Theorems \ref{th:local_asymp_curstat_MSLE} to \ref{th:SLR_test_same_g} be satisfied.
Then
\begin{align}
\label{V_N-representation}
V_n&=\a_N\b_N\int_{t\in[a,b]}\frac{\bigl\{\tilde g_{N2}(t)\tilde h_{N1}(t)-\tilde g_{N1}(t)\tilde h_{N2}(t)\bigr\}^2}{F(t)\{1-F(t)\}\bar g_N(t)g_1(t)g_2(t)}\,dt+O_p\left(N^{-3(1-\a)/2}(\log N)^{3/2}\right),
\end{align}
where
$$
\bar g_N=\a_N g_1+\b_Ng_2.
$$
Moreover,
\begin{align}
\label{V_N-representation2}
&\a_N\b_N\int_{t\in[a,b]}\frac{\bigl\{\tilde g_{N2}(t)\tilde h_{N1}(t)-\tilde g_{N1}(t)\tilde h_{N2}(t)\bigr\}^2}{F(t)\{1-F(t)\}\bar g_N(t)g_1(t)g_2(t)}\,dt-\frac{b-a}{Nb_N}\int K(u)^2\,du\nonumber\\
&=A_N+B_N-C_n+D_N+o_p\left(\frac1{N\sqrt{b_N}}\right),
\end{align}
where
\begin{align*}
A_N
=\frac{2\a_N\b_N}{m^2}\sum_{1<i<j\le m}\bigl\{\dd_i-F(T_i)\bigr\}\bigl\{\dd_j-F(T_j)\bigr\}
\int_{t=a}^b\frac{g_2(t)K_{b_N}(t-T_i)K_{b_N}(t-T_j)}{g_1(t)\bar g_N(t)F(t)\{1-F(t)\}}\,dt,
\end{align*}
\begin{align*}
B_N=\frac{2\a_N\b_N}{n^2}\sum_{m<i<j\le N}\bigl\{\dd_i-F(T_i)\bigr\}\bigl\{\dd_j-F(T_j)\bigr\}
\int_{t=a}^b\frac{g_1(t)K_{b_N}(t-T_i)K_{b_N}(t-T_j)}{g_2(t)\bar g_N(t)F(t)\{1-F(t)\}}\,dt,
\end{align*}
\begin{align*}
C_N=\frac{2\a_N\b_N}{mn}\sum_{i=1}^m\sum_{j=m+1}^N\bigl\{\dd_i-F(T_i)\bigr\}\bigl\{\dd_j-F(T_j)\bigr\}\int_{t=a}^b \frac{K_{b_N}(t-T_i)K_{b_N}(t-T_j)}{\bar g_N(t)F(t)\{1-F(t)\}}\,dt,
\end{align*}
and
\begin{align*}
D_N=\a_N\b_N\int_{t=a}^b\frac{f(t)^2\left\{g_1'(t)g_2(t)-g_2'(t)g_1(t)\right\}^2}{F(t)\{1-F(t)\}\bar g_N(t)g_1(t)g_2(t)}\,dt\left\{\int u^2 K(u)\,du\right\}^2b_N^4.
\end{align*}
Note that $D_N=0$ if $g_1=g_2$.
\end{lemma}
 
\vspace{0.3cm}
\noindent
\vspace{0.3cm}
\noindent
{\bf Proof}.
By Lemma \ref{Talagrand_lemma} and an expansion of the logarithm we get:
\begin{align*}
&2\int_{t\in[a,b]} \left\{\tilde h_{N1}(t)\log \frac{{\tilde F}_{N1}(t)}{{\tilde F}_N(t)}+\bigl\{\tilde g_{N1}(t)-\tilde h_{N1}(t)\bigr\}\log\frac{1- {\tilde F}_{N1}(t)}{1- {\tilde F}_N(t)}\right\}\,dt\\
&=-2\int_{t\in[a,b]} \left\{\tilde h_{N1}(t)\log \frac{{\tilde F}_N(t)}{{\tilde F}_{N1}(t)}+\bigl\{\tilde g_{N1}(t)-\tilde h_{N1}(t)\bigr\}\log\frac{1- {\tilde F}_N(t)}{1- {\tilde F}_{N1}(t)}\right\}\,dt\\
&=\int_{t\in[a,b]}\frac{\tilde g_{N1}(t)^2\bigl\{\tilde F_{N1}(t)-\tilde F_N(t)\bigr\}^2}{\tilde h_{N1}(t)\bigl\{\tilde g_{N1}(t)-\tilde h_{N1}(t)\bigr\}}\,\tilde g_{N1}(t)\,dt+O_p\left(N^{-3(1-\a)/2}(\log N)^{3/2}\right).
\end{align*}
We likewise get, with probability tending to one,
\begin{align*}
&2\int_{t\in[a,b]} \left\{\tilde h_{N2}(t)\log \frac{{\tilde F}_{N2}(t)}{{\tilde F}_N(t)}+\bigl\{\tilde g_{N2}(t)-\tilde h_{N2}(t)\bigr\}\log\frac{1- {\tilde F}_{N2}(t)}{1- {\tilde F}_N(t)}\right\}\,dt\\
&=\int_{t\in[a,b]}\frac{\tilde g_{N2}(t)^2\bigl\{F_{N2}(t)-\tilde F_N(t)\bigr\}^2}{\tilde h_{N2}(t)\bigl\{\tilde g_{N2}(t)-\tilde h_{N2}(t)\bigr\}}\,\tilde g_{N2}(t)\,dt+O_p\left(N^{-3(1-\a)/2}(\log N)^{3/2}\right).
\end{align*}

So we have to consider
\begin{align*}
&\a_N\int_{t\in[a,b]}\frac{\tilde g_{N1}(t)^2\bigl\{\tilde F_{N1}(t)-\tilde F_N(t)\bigr\}^2}{\tilde h_{N1}(t)\bigl\{\tilde g_{N1}(t)-\tilde h_{N1}(t)\bigr\}}\,\tilde g_{N1}(t)\,dt\\
&\qquad\qquad\qquad\qquad+\b_N\int_{t\in[a,b]}\frac{\tilde g_{N2}(t)^2\bigl\{F_{N2}(t)-\tilde F_N(t)\bigr\}^2}{\tilde h_{N2}(t)\bigl\{\tilde g_{N2}(t)-\tilde h_{N2}(t)\bigr\}}\,\tilde g_{N2}(t)\,dt\\
&=\a_N\b_N^2\int_{t\in[a,b]}\frac{\bigl\{\tilde g_{N2}(t)\tilde h_{N1}(t)-\tilde g_{N1}(t)\tilde h_{N2}(t)\bigr\}^2}{\tilde h_{N1}(t)\bigl\{\tilde g_{N1}(t)-\tilde h_{N1}(t)\bigr\}\tilde g_N(t)^2}\,\tilde g_{N1}(t)\,dt\\
&\qquad\qquad\qquad\qquad+\a_N^2\b_N\int_{t\in[a,b]}\frac{\bigl\{\tilde g_{N2}(t)\tilde h_{N1}(t)-\tilde g_{N1}(t)\tilde h_{N2}(t)\bigr\}^2}{\tilde h_{N2}(t)\bigl\{\tilde g_{N2}(t)-\tilde h_{N2}(t)\bigr\}\tilde g_N(t)^2}\,\tilde g_{N2}(t)\,dt
\end{align*}

We have:
\begin{align*}
&\frac{\b_N\tilde g_{N1}(t)}{\tilde h_{N1}(t)\bigl\{\tilde g_{N1}(t)-\tilde h_{N1}(t)\bigr\}\tilde g_N(t)^2}+\frac{\a_N\tilde g_{N2}(t)}{\tilde h_{N2}(t)\bigl\{\tilde g_{N2}(t)-\tilde h_{N2}(t)\bigr\}\tilde g_N(t)^2}\\
&=\frac{\b_N\tilde g_{N1}(t)\tilde h_{N2}(t)\bigl\{\tilde g_{N2}(t)-\tilde h_{N2}(t)\bigr\}+\a_N\tilde g_{N2}(t)\tilde h_{N1}(t)\bigl\{\tilde g_{N1}(t)-\tilde h_{N1}(t)\bigr\}}{\tilde h_{N1}(t)\bigl\{\tilde g_{N1}(t)-\tilde h_{N1}(t)\bigr\}\tilde h_{N2}(t)\bigl\{\tilde g_{N2}(t)-\tilde h_{N2}(t)\bigr\}\tilde g_N(t)^2}\\
&=\frac1{F(t)\{1-F(t)\}\bar g_N(t)g_1(t)g_2(t)}+O_p\left(N^{-(1-\a)/2}\sqrt{\log N}\right).
\end{align*}
Hence:
\begin{align*}
V_n&=\a_N\b_N\int_{t\in[a,b]}\frac{\bigl\{\tilde g_{N2}(t)\tilde h_{N1}(t)-\tilde g_{N1}(t)\tilde h_{N2}(t)\bigr\}^2}{F(t)\{1-F(t)\}\bar g_N(t)g_1(t)g_2(t)}\,dt+O_p\left(N^{-3(1-\a)/2}(\log N)^{3/2}\right).
\end{align*}
Furthermore,
\begin{align*}
&\tilde g_{N2}(t)\tilde h_{N1}(t)-\tilde g_{N1}(t)\tilde h_{N2}(t)\\
&=n^{-1}\sum_{i=m+1}^N K_{b_N}(t-T_i)\,m^{-1}\sum_{i=1}^m \dd_iK_{b_N}(t-T_i)\\
&\qquad\qquad-m^{-1}\sum_{i=1}^m K_{b_N}(t-T_i)\,n^{-1}\sum_{i=m+1}^n \dd_i K_{b_N}(t-T_i)\\
&=n^{-1}\sum_{i=m+1}^N K_{b_N}(t-T_i)\,m^{-1}\sum_{i=1}^m \{F(T_i)-F(t)\}K_{b_N}(t-T_i)\\
&\qquad\qquad-m^{-1}\sum_{i=1}^m K_{b_N}(t-T_i)\,n^{-1}\sum_{i=m+1}^n \{F(T_i)-F(t)\} K_{b_N}(t-T_i)\\
&\qquad\qquad+n^{-1}\sum_{i=m+1}^N K_{b_N}(t-T_i)\,m^{-1}\sum_{i=1}^m\bigl\{\dd_i-F(T_i)\bigr\}K_{b_N}(t-T_i)\\
&\qquad\qquad\qquad\qquad-m^{-1}\sum_{i=1}^m K_{b_N}(t-T_i)\,n^{-1}\sum_{i=m+1}^n \bigl\{\dd_i-F(T_i)\bigr\}K_{b_N}(t-T_i).
\end{align*}

We first consider the first two terms on the right-hand side:
\begin{align*}
&n^{-1}\sum_{i=m+1}^N K_{b_N}(t-T_i)\,m^{-1}\sum_{i=1}^m \bigl\{F(T_i)-F(t)\bigr\}K_{b_N}(t-T_i)\\
&\qquad\qquad-m^{-1}\sum_{i=1}^m K_{b_N}(t-T_i)\,n^{-1}\sum_{i=m+1}^n \bigl\{F(T_i)-F(t)\bigr\} K_{b_N}(t-T_i)\\
&=\tilde g_{N2}(t)\,m^{-1}\sum_{i=1}^m \bigl\{F(T_i)-F(t)\bigr\}K_{b_N}(t-T_i)\\
&\qquad\qquad-\tilde g_{N1}(t)\,n^{-1}\sum_{i=m+1}^n \bigl\{F(T_i)-F(t)\bigr\} K_{b_N}(t-T_i)\\
&=f(t)\left\{\tilde g_{N2}(t)\,m^{-1}\sum_{i=1}^m \{T_i-t\}K_{b_N}(t-T_i)
-\tilde g_{N1}(t)\,n^{-1}\sum_{i=m+1}^n \{T_i-t\} K_{b_N}(t-T_i)\right\}\\
&\qquad\qquad+\tfrac12\tilde g_{N2}(t)\,m^{-1}\sum_{i=1}^m f'(\theta_i)\{T_i-t\}^2K_{b_N}(t-T_i)\\
&\qquad\qquad\qquad\qquad-\tfrac12\tilde g_{N1}(t)\,n^{-1}\sum_{i=m+1}^n f'(\theta_i)\{T_i-t\}^2 K_{b_N}(t-T_i),
\end{align*}
where $\theta_i$ is a point between $t$ and $T_i$. This implies, using the fact that the variance is of order $O(N^{-1}b_N)$,
\begin{align}
\label{difference_bias}
&n^{-1}\sum_{i=m+1}^N K_{b_N}(t-T_i)\,m^{-1}\sum_{i=1}^m \bigl\{F(T_i)-F(t)\bigr\}K_{b_N}(t-T_i)\nonumber\\
&\qquad\qquad-m^{-1}\sum_{i=1}^m K_{b_N}(t-T_i)\,n^{-1}\sum_{i=m+1}^n\bigl\{F(T_i)-F(t)\bigr\} K_{b_N}(t-T_i)\nonumber\\
&=b_N^2f(t)\left\{g_1(t)g_2'(t)-g_1'(t)g_2(t)\right\}\int u^2K(u)\,du\nonumber\\
&\qquad\qquad+\tfrac14b_N^4f(t)\left\{g_1''(t)g_2'(t)-g_1'(t)g_2''(t)\right\}\left\{\int u^2K(u)\,du\right\}^2\nonumber\\
&\qquad\qquad+\tfrac14b_N^4f'(t)\left\{g_1''(t)g_2(t)-g_1(t)g_2''(t)\right\}\left\{\int u^2K(u)\,du\right\}^2+O_p\left(\sqrt{\frac{b_N\log N}{N}}\right)+o\left(b_N^4\right)\nonumber\\
&=b_N^2f(t)\left\{g_1(t)g_2'(t)-g_1'(t)g_2(t)\right\}\int u^2K(u)\,du+O_p\left(\sqrt{\frac{b_N\log N}{N}}\right)+O\left(b_N^4\right),
\end{align}
uniformly for $t\in[a,b]$.

We now define
\begin{equation}
\label{bias}
S_N(t)=b_N^2f(t)\left\{g_1(t)g_2'(t)-g_1'(t)g_2(t)\right\}\int u^2K(u)\,du.
\end{equation}
and
\begin{align}
\label{def_W_N}
&W_N(t)=\tilde g_{N2}(t)\,m^{-1}\sum_{i=1}^m\bigl\{\dd_i-F(T_i)\bigr\}K_{b_N}(t-T_i)\nonumber\\
&\qquad\qquad\qquad\qquad-\tilde g_{N1}(t)\,n^{-1}\sum_{i=m+1}^n \bigl\{\dd_i-F(T_i)\bigr\}
K_{b_N}(t-T_i).
\end{align}
Then
$$
E\left(W_N(t)\bigm|T_1,\dots,T_N\right)=0,\qquad \mbox{var}\left(W_N(t)\right)=O\left(\frac1{Nb_N}\right),
$$
and hence:
\begin{align*}
&\left\{\tilde g_{N2}(t)\tilde h_{N1}(t)-\tilde g_{N1}(t)\tilde h_{N2}(t)\right\}^2
=\left\{W_N(t)+S_N(t)\right\}^2+O_p\left(\frac{\log N}{N}\right)+O_p\left(b_N^6\right).
\end{align*}
We have:
\begin{align*}
\a_N\b_N\int_{t\in[a,b]}\frac{S_N(t)W_N(t)}{F(t)\{1-F(t)\}\bar g_N(t)g_1(t)g_2(t)}\,dt
=O_p\left(\frac{b_N^2}{\sqrt{N}}\right),
\end{align*}
since, by the central limit theorem
$$
\int_{t\in[a,b]}\frac{S_N(t)W_N(t)}{F(t)\{1-F(t)\}\bar g_N(t)g_1(t)g_2(t)}\,dt=O_p\left(N^{-1/2}\right).
$$
Note that this term is zero if $g_1=g_2$.

So we get:
\begin{align}
\label{V_N-expansion}
&\a_N\b_N\int_{t\in[a,b]}\frac{\bigl\{\tilde g_{N2}(t)\tilde h_{N1}(t)-\tilde g_{N1}(t)\tilde h_{N2}(t)\bigr\}^2}{F(t)\{1-F(t)\}\bar g_N(t)g_1(t)g_2(t)}\,dt\nonumber\\
&=\a_N\b_N\int_{t\in[a,b]}\frac{W_N(t)^2}{F(t)\{1-F(t)\}\bar g_N(t)g_1(t)g_2(t)}\,dt+D_N\nonumber\\
&\qquad\qquad\qquad\qquad+O_p\left(\frac{b_N^2}{\sqrt{N}}\right)+O_p\left(\frac{\log N}{N}\right)+O_p\left(b_N^6\right),
\end{align}
where $D_N$ is defined as in the formulation of the lemma, and where the term $O_p\left(b_N^2/\sqrt{N}\right)$ is absent if $g_1=g_2$. Let
\begin{align*}
&W_{N1}(t)=\tilde g_{N2}(t)\,m^{-1}\sum_{i=1}^m\bigl\{\dd_i-F(T_i)\bigr\}K_{b_N}(t-T_i),
\end{align*}
and
\begin{align*}
&W_{N2}(t)=\tilde g_{N1}(t)\,n^{-1}\sum_{i=m+1}^n \bigl\{\dd_i-F(T_i)\bigr\}
K_{b_N}(t-T_i).
\end{align*}
Then, by definition (\ref{def_W_N}), $W_N=W_{N1}+W_{N2}$, and we get:
\begin{align*}
&\a_N\b_N\int_{t\in[a,b]}\frac{W_N(t)^2}{F(t)\{1-F(t)\}\bar g_N(t)g_1(t)g_2(t)}\,dt\\
&=\a_N\b_N\int_{t\in[a,b]}\frac{W_{N1}(t)^2}{F(t)\{1-F(t)\}\bar g_N(t)g_1(t)g_2(t)}\,dt\\
&\qquad+\a_N\b_N\int_{t\in[a,b]}\frac{W_{N2}(t)^2}{F(t)\{1-F(t)\}\bar g_N(t)g_1(t)g_2(t)}\,dt\\
&\qquad\qquad\qquad-2\a_N\b_N\int_{t\in[a,b]}\frac{W_{N1}(t)W_{N2}(t)}{F(t)\{1-F(t)\}\bar g_N(t)g_1(t)g_2(t)}\,dt
\end{align*}
We now have, using Lemma \ref{Talagrand_lemma} for $\tilde g_{N2}$,
\begin{align*}
&\a_N\b_N\int_{t\in[a,b]}\frac{W_{N1}(t)^2}{F(t)\{1-F(t)\}\bar g_N(t)g_1(t)g_2(t)}\,dt\\
&=\a_N\b_N\int_{t\in[a,b]}\frac{g_2(t)\left\{m^{-1}\sum_{i=1}^m K_{b_N}(t-T_i)\left\{\dd_i-F(T_i)\right\}\right\}^2}{F(t)\{1-F(t)\}\bar g_N(t)g_1(t)}\,dt+o_p\left(\frac1{N\sqrt{b_N}}\right)\\
&=\frac{\a_N\b_N}{m^2}\sum_{i=1}^m \bigl\{\dd_i-F(T_i)\bigr\}^2\int_{t=a}^b \frac{g_2(t)K_{b_N}(t-T_i)^2}{g_1(t)\bar g_N(t)F(t)\{1-F(t)\}}\,dt\\
&\qquad+\frac{2\a_N\b_N}{m^2}\sum_{1<i<j\le m}\bigl\{\dd_i-F(T_i)\bigr\}\bigl\{\dd_j-F(T_j)\bigr\}
\int_{t=a}^b\frac{g_2(t)K_{b_N}(t-T_i)K_{b_N}(t-T_j)}{g_1(t)\bar g_N(t)F(t)\{1-F(t)\}}\,dt\\
&\qquad\qquad\qquad\qquad\qquad\qquad\qquad\qquad\qquad\qquad\qquad\qquad\qquad\qquad\qquad+o_p\left(\frac1{N\sqrt{b_N}}\right).
\end{align*}
Moreover, by the central limit theorem,
\begin{align*}
&\frac{\a_N\b_N}{m^2}\sum_{i=1}^m \bigl\{\dd_i-F(T_i)\bigr\}^2\int_{t=a}^b \frac{g_2(t)K_{b_N}(t-T_i)^2}{g_1(t)\bar g_N(t)F(t)\{1-F(t)\}}\,dt\\
&\qquad\qquad\qquad\qquad-\frac{\a_N\b_N}{m^2}E\sum_{i=1}^m \bigl\{\dd_i-F(T_i)\bigr\}^2\int_{t=a}^b \frac{g_2(t)K_{b_N}(t-T_i)^2}{g_1(t)\bar g_N(t)F(t)\{1-F(t)\}}\,dt\\
&=O_p\left(\frac1{N^{3/2}b_N}\right)=o_p\left(\frac1{N\sqrt{b_N}}\right),
\end{align*}
and
\begin{align*}
&\frac{\a_N\b_N}{m^2}E\sum_{i=1}^m E\bigl\{\dd_i-F(T_i)\bigr\}^2\int_{t=a}^b \frac{g_2(t)K_{b_N}(t-T_i)^2}{g_1(t)\bar g_N(t)F(t)\{1-F(t)\}}\,dt\\
&=\frac{\b_N}{N m}\sum_{i=1}^m E\left\{\bigl\{\dd_i-F(T_i)\bigr\}^2\int_{t=a}^b \frac{g_2(t)K_{b_N}(t-T_i)^2}{g_1(t)\bar g_N(t)F(t)\{1-F(t)\}}\,dt\right\}\\
&=\frac{\b_N}{N}E\int_{t=a}^b \frac{g_2(t)K_{b_N}(t-T_1)^2}{g_1(t)\bar g_N(t)}\,dt
=\frac1{Nb_N}\int_{t=a}^b \frac{\b_Ng_2(t)}{\bar g_N(t)}\,dt\int K(u)^2\,du+O\left(\frac{b_N}{N}\right).
\end{align*}
We similarly get:
\begin{align*}
&\frac{\a_N\b_N}{n^2}\sum_{i=m+1}^N \bigl\{\dd_i-F(T_i)\bigr\}^2\int_{t=a}^b \frac{g_1(t)K_{b_N}(t-T_i)^2}{g_2(t)\bar g_N(t)F(t)\{1-F(t)\}}\,dt\\
&\qquad\qquad\qquad\qquad-\frac{\a_N\b_N}{n^2}E\sum_{i=m+1}^N \bigl\{\dd_i-F(T_i)\bigr\}^2\int_{t=a}^b \frac{g_1(t)K_{b_N}(t-T_i)^2}{g_2(t)\bar g_N(t)F(t)\{1-F(t)\}}\,dt\\
&=O_p\left(\frac1{N^{3/2}b_N}\right)=o_p\left(\frac1{N\sqrt{b_N}}\right),
\end{align*}
and
\begin{align*}
&\frac{\a_N\b_N}{m^2}E\sum_{i=m+1}^N \bigl\{\dd_i-F(T_i)\bigr\}^2\int_{t=a}^b \frac{g_2(t)K_{b_N}(t-T_i)^2}{g_1(t)\bar g_N(t)F(t)\{1-F(t)\}}\,dt\\
&=\frac1{Nb_N}\int_{t=a}^b \frac{\a_Ng_1(t)}{\bar g_N(t)}\,dt\int K(u)^2\,du+O\left(\frac{b_N}{N}\right).
\end{align*}
Hence:
\begin{align*}
&\frac{\a_N\b_N}{m^2}\sum_{i=1}^m \bigl\{\dd_i-F(T_i)\bigr\}^2\int_{t=a}^b \frac{g_2(t)K_{b_N}(t-T_i)^2}{g_1(t)\bar g_N(t)F(t)\{1-F(t)\}}\,dt\\
&\qquad\qquad+\frac{\a_N\b_N}{n^2}\sum_{i=m+1}^N \bigl\{\dd_i-F(T_i)\bigr\}^2\int_{t=a}^b \frac{g_1(t)K_{b_N}(t-T_i)^2}{g_2(t)\bar g_N(t)F(t)\{1-F(t)\}}\,dt\\
&=\frac1{Nb_N}\int_{t=a}^b \frac{\a_Ng_1(t)+\b_Ng_2(t)}{\bar g_N(t)}\,dt\int K(u)^2\,du+o_p\left(\frac1{N\sqrt{b_N}}\right)\\
&=\frac{b-a}{Nb_N}\int K(u)^2\,du+o_p\left(\frac1{N\sqrt{b_N}}\right).
\end{align*}
The representation (\ref{V_N-representation2}) now follows.\eop

\vspace{0.3cm}
\noindent
{\bf Proofs of Theorem \ref{th:local_asymp_curstat_MSLE} to \ref{th:SLR_test_same_g}}. By Lemma \ref{expectation_lemma}, we only have to study the terms on the right-hand side of (\ref{V_N-representation2}). We condition on the values $T_1,\dots,T_N$. The  first term $A_N$ can be written
\begin{align*}
A_N
=\sum_{j=2}^m Y_j,
\end{align*}
where
\begin{align*}
Y_j=\frac{2\a_N\b_N}{m^2}\sum_{1\le i<j}\bigl\{\dd_i-F(T_i)\bigr\}\bigl\{\dd_j-F(T_j)\bigr\}
\int_{t=a}^b\frac{g_2(t)K_{b_N}(t-T_i)K_{b_N}(t-T_j)}{g_1(t)\bar g_N(t)F(t)\{1-F(t)\}}\,dt.
\end{align*}
Letting ${\cal F}_j$ be the $\s$-algebra, generated by $Y_1,\dots,Y_j$, and ${\cal F}_0$ be the trivial $\s$-algebra, we get:
$$
E\left\{Y_j\bigm| {\cal F}_{j-1}\right\}=0,\,j=1,\dots,m.
$$
Furthermore we have, in probability,
\begin{align*}
&E\left\{Y_j^2\bigm| {\cal F}_{j-1}\right\}\\
&= \frac{4\a_N^2\b_N^2}{m^4}F(T_j)\bigl\{1-F(T_j)\bigr\}\left\{\sum_{1\le i<j}\bigl\{\dd_i-F(T_i)\bigr\}
\int_{t=a}^b\frac{g_2(t)K_{b_N}(t-T_i)K_{b_N}(t-T_j)}{g_1(t)\bar g_N(t)F(t)\{1-F(t)\}}\,dt\right\}^2\\
&\sim \frac{4(j-1)\a_N^2\b_N^2}{m^4b_N}\frac{g_2(T_j)^21_{[a-b_N,b+b_N]}(T_j)}{\bar g_N(T_j)^2g_1(T_j)}\int\left\{\int K(v)K(v+x)\,dv\right\}^2\,dx,
\end{align*}
where the last relation holds for large $j$.
Hence we get, in probability,
\begin{align*}
&\sum_{j=1}^m E\left\{Y_j^2\bigm| {\cal F}_{j-1}\right\}
\sim\sum_{j=1}^m\frac{4(j-1)\b_N^2}{m^2N^2b_N}\frac{g_2(T_j)^21_{[a-b_N,b+b_N]}(T_j)}{\bar g_N(T_j)^2g_1(T_j)}\int\left\{\int K(v)K(v+x)\,dv\right\}^2\,dx\\
&\sim \frac{2m(m-1)\b_N^2}{m^2N^2b_N}\int_a^b\frac{g_2(t)^2}{\bar g_N(t)^2}\,dt\int\left\{\int K(v)K(v+x)\,dv\right\}^2\,dx\\
&\sim \frac{2\b_N^2}{N^2b_N}\int_a^b\frac{g_2(t)^2}{\bar g_N(t)^2}\,dt\int\left\{\int K(v)K(v+x)\,dv\right\}^2\,dx,\,m\to\infty.
\end{align*}
We use here that (for the $T_j$ being random again):
\begin{align*}
&\mbox{var}\left(\sum_{j=1}^m\frac{4(j-1)}{m^2}\frac{g_2(T_j)^21_{[a-b_N,b+b_N]}(T_j)}{\bar g_N(T_j)^2g_1(T_j)}\right)\\
&=\frac1{m^4}\mbox{var}\left(\frac{4g_2(T_1)^21_{[a-b_N,b+b_N]}(T_1)}{\bar g_N(T_1)^2g_1(T_1)}\right)\sum_{j=1}^m (j-1)^2=O\left(\frac1N\right),
\end{align*}
implying
\begin{align*}
&\sum_{j=1}^m\frac{4(j-1)}{m^2}\frac{g_2(T_j)^21_{[a-b_N,b+b_N]}(T_j)}{\bar g_N(T_j)^2g_1(T_j)}
\sim \sum_{j=1}^m\frac{4(j-1)}{m^2}E\left\{\frac{g_2(T_j)^21_{[a-b_N,b+b_N]}(T_j)}{\bar g_N(T_j)^2g_1(T_j)}\right\}\\
&\sim \sum_{j=1}^m\frac{4(j-1)}{m^2}\int_a^b\frac{g_2(t)^2}{\bar g_N(t)^2}\,dt
\sim \frac{2m(m-1)}{m^2}\int_a^b\frac{g_2(t)^2}{\bar g_N(t)^2}\,dt\sim 2\int_a^b\frac{g_2(t)^2}{\bar g_N(t)^2}\,dt.
\end{align*}
By similar methods we can extend this to the indices $j=m+1,\dots,N$, where
\begin{align*}
Y_{j}=&\frac{2\a_N\b_N}{n^2}\sum_{m+1\le i<j}\bigl\{\dd_i-F(T_i)\bigr\}\bigl\{\dd_j-F(T_j)\bigr\}
\int_{t=a}^b\frac{g_2(t)K_{b_N}(t-T_i)K_{b_N}(t-T_j)}{g_1(t)\bar g_N(t)F(t)\{1-F(t)\}}\,dt\\
&+\frac{2\a_N\b_N}{mn}\sum_{i=1}^m\bigl\{\dd_i-F(T_i)\bigr\}\bigl\{\dd_j-F(T_j)\bigr\}
\int_{t=a}^b\frac{K_{b_N}(t-T_i)K_{b_N}(t-T_j)}{\bar g_N(t)F(t)\{1-F(t)\}}\,dt,
\end{align*}
which also involves the terms $B_N$ and $C_N$, and results in:
\begin{align*}
&\sum_{j=1}^N E\left\{Y_j^2\bigm| {\cal F}_{j-1}\right\}\\
&\sim\frac2{N^2b_N}\left\{\int  K(v)K(v+x)\,dv\right\}^2\,dx
\int_{t=a}^b\frac{\b_N^2g_2(t)^2+\a_N^2g_1(t)^2+2\a_N\b_N g_1(t)g_2(t)}{\bar g_N(t)^2}\,dt\\
&=\frac{2(b-a)}{N^2b_N}\left\{\int  K(v)K(v+x)\,dv\right\}^2\,dx.
\end{align*}
So we find:
\begin{align*}
&\sum_{j=1}^N E\left\{\left(N\sqrt{b_N}Y_j\right)^2\bigm| {\cal F}_{j-1}\right\}
\stackrel{p}\longrightarrow 2(b-a)\left\{\int  K(v)K(v+x)\,dv\right\}^2\,dx,\,N\to\infty.
\end{align*}
By tedious but straightforward computations, using $4$th moments of the Bernoulli distribution, one can also check that
\begin{align*}
&\sum_{j=1}^N E\left\{\left(N\sqrt{b_N}Y_j\right)^21_{\{N^2b_NY_j^2>\e\}}\bigm| {\cal F}_{j-1}\right\}
\stackrel{p}\longrightarrow 0,\,N\to\infty.
\end{align*}
The result now follows from the martingale convergence theorem on p.\ 171 in \shortciteN{pollard:84}.\eop

\vspace{0.3cm}
\noindent
{\bf Sketch of proof of (\ref{Andersen_asymp}).}
 First consider
$$
\int_0^a\left\{\hat F_m(t)^2-F(t)^2\right\}\,dG(t),
$$
where we assume $G_1=G_2=G$. Then:
\begin{align*}
\int_0^a\left\{\hat F_m(t)^2-F(t)^2\right\}\,dG(t)
&=2\int_0^a\left\{\hat F_m(t)-F(t)\right\}F(t)\,dG(t)+\int_0^a\left\{\hat F_m(t)-F(t)\right\}^2\,dG(t)\\
&=2\int_0^a\left\{\hat F_m(t)-F(t)\right\}F(t)\,dG(t)+O_p\left(m^{-2/3}\right).
\end{align*}
Secondly,
\begin{align*}
&2\int_0^a\left\{\hat F_m(t)-F(t)\right\}F(t)\,dG(t)
=2\int_0^a\left\{\hat F_m(t)-\d\right\}F(t)\,dP_{01}(t,\d),
\end{align*}
where $P_{01}$ is the probability measure, generating the random variables $(T_1,\dd_1),\dots,(T_m,\dd_m)$. Let $\bar F$ be a piecewise constant version of $F$, which is constant on the same intervals as $\hat F_m$. Then:
\begin{align*}
&2\int_0^a\bigl\{\hat F_m(t)-\d\bigr\}F(t)\,dP_{01}(t,\d)\\
&=2\int_0^a\bigl\{\hat F_m(t)-\d\bigr\}\bar{F}_0(t)\,dP_{01}(t,\d)
+2\int_0^a\bigl\{\hat F_m(t)-\d\bigr\}\bigl\{F(t)-\bar{F}_0(t)\bigr\}\,dP_{01}(t,\d)\\
&=2\int_0^a\bigl\{\hat F_m(t)-\d\bigr\}\bar{F}_0(t)\,dP_{01}(t,\d)
+2\int_0^a\bigl\{\hat F_m(t)-F(t)\bigr\}\bigl\{F(t)-\bar{F}_0(t)\bigr\}\,dG(t)\\
&=2\int_0^a\bigl\{\hat F_m(t)-\d\bigr\}\bar{F}_0(t)\,dP_{01}(t,\d)+O_p\left(m^{-2/3}\right).
\end{align*}
But, by the characterization of the MLE $\hat F_m$, we have, if $\t(a)$ is the last point of jump of $\hat F_m$ before $a$,
$$
2\int_{[0,\t(a))}\bigl\{\hat F_m(t)-\d\bigr\}\bar{F}_0(t)\,d\P_{N1}(t,\d)=0,
$$
and hence:
\begin{align*}
&2\int_0^a\bigl\{\hat F_m(t)-\d\bigr\}\bar{F}_0(t)\,dP_{01}(t,\d)
=2\int_{[0,\t(a))}\bigl\{\hat F_m(t)-\d\bigr\}\bar{F}_0(t)\,d\left(P_{01}-\P_{N1}\right)(t,\d)\\
&\qquad\qquad\qquad\qquad\qquad\qquad\qquad\qquad\qquad\qquad\qquad\qquad\qquad\qquad+O_p\left(m^{-2/3}\right)\\
&=2\int_{[0,\t(a))}\bigl\{F(t)-\d\bigr\}\bar{F}_0(t)\,d\left(P_{01}-\P_{N1}\right)(t,\d)\\
&\qquad\qquad\qquad\qquad+2\int_{[0,\t(a))}\bigl\{\hat F_m(t)-F(t)\bigr\}\bar{F}_0(t)\,d\left(P_{01}-\P_{N1}\right)(t,\d)+O_p\left(m^{-2/3}\right)\\
&=2\int_{[0,a]}\bigl\{F(t)-\d\bigr\}F(t)\,d\left(P_{01}-\P_{N1}\right)(t,\d)+O_p\left(m^{-2/3}\right),
\end{align*}
where the first term, multiplied by $\sqrt{m}$, is asymptotically normal with mean zero and variance
$$
4\int_0^a F(t)^3\bigl\{1-F(t)\bigr\}\,dG(t).
$$
This implies the result, since we can write:
\begin{align*}
&\int_0^a\left\{\hat F_m(t)^2-\hat F_n(t)^2\right\}\,d\G_N(t)\\
&=\int_0^a\left\{\hat F_m(t)^2-\hat F_n(t)^2\right\}\,dG(t)
+\int_{[0,a]}\left\{\hat F_m(t)^2-\hat F_n(t)^2\right\}\,d\left(\G_N-G\right)(t)\\
&=\int_0^a\left\{\hat F_m(t)^2-F(t)^2\right\}\,dG(t)
-\int_0^a\left\{\hat F_n(t)^2-F(t)^2\right\}\,dG(t)+O_p\left(N^{-2/3}\right),
\end{align*}
and since $\hat F_m$ and $\hat F_n$ are based on independent samples.\eop

\vspace{0.3cm}
\noindent
{\bf Proof of Theorem \ref{th:bootstrap1}}. We may assume that, for large $N$, $\tilde F_{N,\tilde b_N}$ has the representation
$$
\tilde F_{N,\tilde b_N}(t)=\frac{\int \d K_{\tilde b_N}(t-u)\,d\P_N(u,\d)}{\int K_{\tilde b_N}(t-u)\,d\G_N(u)}\,.
$$
for $t\in[a,b]$, where $\tilde b_N\asymp N^{-1/5}$. This gives
$$
\tilde f_{N,\tilde b_N}(t)=\frac{\int \d K_{\tilde b_N}'(t-u)\,d\P_N(u,\d)}{\tilde g_{N,\tilde b_N}(t)}
-\frac{\tilde g_{N,\tilde b_N}'(t)\int \d K_{\tilde b_N}(t-u)\,d\P_N(u,\d)}{\tilde g_{N,\tilde b_N}(t)^2}\,,
$$
where
$$
\tilde g_{N,\tilde b_N}(t)=\int K_{\tilde b_N}(t-u)\,d\G_N(u),\qquad\tilde g_{N,\tilde b_N}'(t)=\int K_{\tilde b_N}'(t-u)\,d\G_N(u),
$$
and
$$
K_{\tilde b_N}'(t-u)=\frac1{\tilde b_N^2}K'\left(\frac{t-u}{\tilde b_N}\right)
$$
By the assumptions on $g$, and using $\tilde b_N\asymp n^{-1/5}$, we have
\begin{align*}
&\sup_{t\in[a,b]}\left|\tilde g_{N,\tilde b_N}(t)-\bar g_N(t)\right|=O_p\left(N^{-2/5}\sqrt{\log n}\right)
\mbox{ and }\sup_{t\in[a,b]}\left|\tilde g_{N,\tilde b_N}'(t)-g'(t)\right|\\
&=O_p\left(N^{-1/5}\sqrt{\log n}\right)
\end{align*}
uniformly for $t\in[a,b]$. Furthermore, since
$$
\int \d K_{\tilde b_N}'(t-u)\,d\P_N(u,\d)=\frac1{N\tilde b_N^2}\sum_{i=1}^N K'\left(\frac{t-T_i}{\tilde b_N}\right)\dd_i
$$
we get:
\begin{align*}
&\int \d K_{\tilde b_N}'(t-u)\,d\P_N(u,\d)-\int K_{\tilde b_N}'(t-u)F(u)\,dG(u)\\
&=\int \left\{\d-F(u)\right\} K_{\tilde b_N}'(t-u)\,d\P_N(u,\d)
+\int F(u) K_{\tilde b_N}'(t-u)\,d\left(\G_N-G\right)(u)
\end{align*}
and hence
\begin{equation}
\label{uniform_conv_derivative}
\sup_{t\in[a,b]}\left|\tilde f_{N,\tilde b_N}(t)-f(t)\right|=O_p\left(N^{-1/5}\sqrt{\log N}\right).
\end{equation}
It can be proved in a similar way that
$$
\sup_{x\in[a,b]}\left|\tilde F_{N,\tilde b_N}(t)-F(t)\right|=O_p\left(N^{-2/5}\sqrt{\log N}\right).
$$

The bootstrap test statistic $V_N^*$ now has the representation
\begin{align*}
&V_N=\frac{2m}{N}\int_{t\in[a,b]} \left\{\tilde h_{N1}^*(t)\log \frac{{\tilde F}_{N1}^*(t)}{{\tilde F}_N^*(t)}+\bigl\{\tilde g_{N1}(t)-\tilde h_{N1}^*(t)\bigr\}\log\frac{1- {\tilde F}_{N1}^*(t)}{1- {\tilde F}_N^*(t)}\right\}\,dt\nonumber\\
&\qquad\qquad+\frac{2n}{N}\int_{t\in[a,b]} \left\{\tilde h_{N2}^*(t)\log \frac{{\tilde F}_{N2}^*(t)}{\tilde F_N^*(t)}+\bigl\{\tilde g_{N2}(t)-\tilde h_{N2}^*(t)\bigr\}\log\frac{1- {\tilde F}_{N2}^*(T_i)}{1- {\tilde F}_N^*(T_i)}\right\}\,dt,
\end{align*}
where
\begin{align*}
\tilde h_{Nj}^*(t)=\int \d^*K_{b_N}(t-u)\,d\P_{Nj}(u,\d^*),\,j=1,2,
\end{align*}
and the $\dd_i^*$ are defined by
$$
\dd_i^*=1_{\left[0,\tilde F_{N,\tilde b_N}(T_i)\right]}(U_i^*),
$$
for independent random variables $U_1^*,\dots,U_N^*$, independent of the random variables $(T_i,\dd_i)$, $i=1,\dots,N$, and where we may assume, as before, that
$$
{\tilde F}_{Nj}^*(t)
=\frac{\int \d^*K_{b_N}(t-u)\,d\P_{Nj}(u,\d^*)}{\tilde g_{Nj}(t)},\,j=1,2.
$$
Note that the only extra randomness is introduced by the uniform random variables $U_i^*$, and that the bandwidth $b_N$, used here, may be smaller than the bandwidth $\tilde b_N$, used in the computation of $\tilde F_{N,\tilde b_N}$. In fact, $b_N$ is the bandwidth which is used in the original sample and we have, by assumption
$$
b_N\asymp N^{-\a},
$$
where $1/3<\a<1/5$, and where we allow $\a=1/5$ if it is assumed that $g_1=g_2$. The densities $\tilde g_{Nj}$ have been computed in the original sample, using this possibly smaller bandwidth $b_N$.

We now get, similarly as in Lemma \ref{expectation_lemma},
\begin{align*}
V_n^*&=\a_N\b_N\int_{t\in[a,b]}\frac{\bigl\{\tilde g_{N2}(t)\tilde h_{N1}^*(t)-\tilde g_{N1}(t)\tilde h_{N2}^*(t)\bigr\}^2}{\tilde F_{N,\tilde b_N}(t)\{1-\tilde F_{N,\tilde b_N}(t)\}\bar g_N(t)g_1(t)g_2(t)}\,dt+o_p\left(\frac1{N\sqrt{b_N}}\right),
\end{align*}
where
$$
\bar g_N=\a_N g_1+\b_Ng_2,
$$
and
\begin{align*}
&\tilde g_{N2}(t)\tilde h_{N1}^*(t)-\tilde g_{N1}(t)\tilde h_{N2}^*(t)\\
&=n^{-1}\sum_{i=m+1}^N K_{b_N}(t-T_i)\,m^{-1}\sum_{i=1}^m \dd_i^*K_{b_N}(t-T_i)\\
&\qquad\qquad-m^{-1}\sum_{i=1}^m K_{b_N}(t-T_i)\,n^{-1}\sum_{i=m+1}^n \dd_i^* K_{b_N}(t-T_i)\\
&=n^{-1}\sum_{i=m+1}^N K_{b_N}(t-T_i)\,m^{-1}\sum_{i=1}^m \bigl\{\tilde F_{N,\tilde b_N}(T_i)-\tilde F_{N,\tilde b_N}(t)\bigr\}K_{b_N}(t-T_i)\\
&\qquad\qquad-m^{-1}\sum_{i=1}^m K_{b_N}(t-T_i)\,n^{-1}\sum_{i=m+1}^n\bigl\{\tilde F_{N,\tilde b_N}(T_i)-\tilde F_{N,\tilde b_N}(t)\bigr\} K_{b_N}(t-T_i)\\
&\qquad\qquad+n^{-1}\sum_{i=m+1}^N K_{b_N}(t-T_i)\,m^{-1}\sum_{i=1}^m\bigl\{\dd_i^*-\tilde F_{N,\tilde b_N}(T_i)\bigr\}K_{b_N}(t-T_i)\\
&\qquad\qquad\qquad\qquad-m^{-1}\sum_{i=1}^m K_{b_N}(t-T_i)\,n^{-1}\sum_{i=m+1}^n \bigl\{\dd_i^*-\tilde F_{N,\tilde b_N}(T_i)\bigr\}K_{b_N}(t-T_i).
\end{align*}
This is the same decomposition as used in the proof of Lemma \ref{expectation_lemma}, but with $F$ replaced by $\tilde F_{N,\tilde b_N}$ and $\dd_i$ replaced by $\dd_i^*$. Instead of $W_N$, defined by (\ref{def_W_N}), we get:
\begin{align*}
&W_N^*(t)=\tilde g_{N2}(t)\,m^{-1}\sum_{i=1}^m\bigl\{\dd_i^*-\tilde F_{N,\tilde b_N}(T_i)\bigr\}K_{b_N}(t-T_i)\nonumber\\
&\qquad\qquad\qquad\qquad-\tilde g_{N1}(t)\,n^{-1}\sum_{i=m+1}^n \bigl\{\dd_i^*-\tilde F_{N,\tilde b_N}(T_i)\bigr\}
K_{b_N}(t-T_i).
\end{align*}
and
\begin{align*}
&\left\{\tilde g_{N2}(t)\tilde h_{N1}^*(t)-\tilde g_{N1}(t)\tilde h_{N2}^*(t)\right\}^2
=\left\{W_N^*(t)+S_N^*(t)\right\}^2+o_p\left(\frac1{N\sqrt{b_N}}\right),
\end{align*}
where
\begin{align*}
S_N^*(t)=b_N^2\tilde f_{N,\tilde b_N}(t)\left\{g_1(t)g_2'(t)-g_1'(t)g_2(t)\right\}\int u^2K(u)\,du.
\end{align*}
Moreover,
\begin{align}
&\a_N\b_N\int_{t\in[a,b]}\frac{\bigl\{\tilde g_{N2}(t)\tilde h_{N1}^*(t)-\tilde g_{N1}(t)\tilde h_{N2}^*(t)\bigr\}^2}{\tilde F_{N,\tilde b_N}(t)\{1-\tilde F_{N,\tilde b_N}(t)\}\bar g_N(t)g_1(t)g_2(t)}\,dt-\frac{b-a}{Nb_N}\int K(u)^2\,du\nonumber\\
&=A_N^*+B_N^*-C_n^*+D_N+o_p\left(\frac1{N\sqrt{b_N}}\right),
\end{align}
where
\begin{align*}
A_N^*
&=\frac{2\a_N\b_N}{m^2}\sum_{1<i<j\le m}\bigl\{\dd_i^*-\tilde F_{N,\tilde b_N}(T_i)\bigr\}\bigl\{\dd_j^*-\tilde F_{N,\tilde b_N}(T_j)\bigr\}\\
&\qquad\qquad\qquad\qquad\qquad\qquad\cdot\int_{t=a}^b\frac{g_2(t)K_{b_N}(t-T_i)K_{b_N}(t-T_j)}{g_1(t)\bar g_N(t)\tilde F_{N,\tilde b_N}(t)\{1-\tilde F_{N,\tilde b_N}(t)\}}\,dt,
\end{align*}
\begin{align*}
B_N^*
&=\frac{2\a_N\b_N}{n^2}\sum_{m<i<j\le N}\bigl\{\dd_i^*-\tilde F_{N,\tilde b_N}(T_i)\bigr\}\bigl\{\dd_j^*-F(T_j)\bigr\}\\
&\qquad\qquad\qquad\qquad\qquad\qquad\cdot
\int_{t=a}^b\frac{g_1(t)K_{b_N}(t-T_i)K_{b_N}(t-T_j)}{g_2(t)\bar g_N(t)\tilde F_{N,\tilde b_N}(t)\{1-\tilde F_{N,\tilde b_N}(t)\}}\,dt,
\end{align*}
\begin{align*}
C_N^*
&=\frac{2\a_N\b_N}{mn}\sum_{i=1}^m\sum_{j=m+1}^N\bigl\{\dd_i^*-\tilde F_{N,\tilde b_N}(T_i)\bigr\}\bigl\{\dd_j^*-\tilde F_{N,\tilde b_N}(T_j)\bigr\}\\
&\qquad\qquad\qquad\qquad\qquad\qquad\cdot\int_{t=a}^b \frac{K_{b_N}(t-T_i)K_{b_N}(t-T_j)}{\bar g_N(t)\tilde F_{N,\tilde b_N}(t)\{1-\tilde F_{N,\tilde b_N}(t)\}}\,dt,
\end{align*}
and the bias term $D_N$ is given by:
\begin{align*}
D_N=\a_N\b_N\int_{t=a}^b\frac{f(t)^2\left\{g_1'(t)g_2(t)-g_2'(t)g_1(t)\right\}^2}{F(t)\{1-F(t)\}\bar g_N(t)g_1(t)g_2(t)}\,dt\left\{\int u^2 K(u)\,du\right\}^2b_N^4.
\end{align*}
Note (again) that $D_N=0$ if $g_1=g_2$.

However, the distribution function $\tilde F_{N,\tilde b_N}$ does not satisfy the condition that the second derivative is uniformly bounded on an interval $(a',b')$, containing $[a,b]$, which is a condition on $F$ in Theorems
\ref{th:local_asymp_curstat_MSLE} to \ref{th:SLR_test_same_g}. But a scrutiny of the proof of Lemma \ref{expectation_lemma} reveals that this condition was only needed to take care of the bias term
\begin{align*}
&n^{-1}\sum_{i=m+1}^N K_{b_N}(t-T_i)\,m^{-1}\sum_{i=1}^m \bigl\{F(T_i)-F(t)\bigr\}K_{b_N}(t-T_i)\nonumber\\
&\qquad\qquad-m^{-1}\sum_{i=1}^m K_{b_N}(t-T_i)\,n^{-1}\sum_{i=m+1}^n \bigl\{F(T_i)-F(t)\bigr\} K_{b_N}(t-T_i),
\end{align*}
see (\ref{difference_bias}), which in the present case transforms into
\begin{align*}
&n^{-1}\sum_{i=m+1}^N K_{b_N}(t-T_i)\,m^{-1}\sum_{i=1}^m \bigl\{\tilde F_{N,\tilde b_N}(T_i)-\tilde F_{N,\tilde b_N}(t)\bigr\}K_{b_N}(t-T_i)\nonumber\\
&\qquad\qquad-m^{-1}\sum_{i=1}^m K_{b_N}(t-T_i)\,n^{-1}\sum_{i=m+1}^n \bigl\{\tilde F_{N,\tilde b_N}(T_i)-\tilde F_{N,\tilde b_N}(t)\bigr\} K_{b_N}(t-T_i),
\end{align*}
since we do not change the $T_i$ of the original samples. But since
$$
\int \d K_{\tilde b_N}''(t-u)\,d\P_N(u,\d)=\frac1{N\tilde b_N^3}\sum_{i=1}^N K''\left(\frac{t-T_i}{\tilde b_N}\right)\dd_i=O_p\left(\sqrt{\log N}\right),
$$
uniformly in $t\in[a,b]$, using again the methods of Lemma \ref{Talagrand_lemma} together with the assumption that $F$, $g_1$ and $g_2$ are twice continuously differentiable, the remainder term $O(b_N^4)$ in (\ref{difference_bias}) can be replaced by a remainder term of order $O_p(b_N^4\log N)$, which is sufficient for our purposes. Theorem \ref{th:bootstrap1} now follows.\eop

\vspace{0.3cm}
\noindent
Piet Groeneboom,
Delft University of Technology,
Department of Applied Mathematics,\\
Mekelweg 4,
2628 CD Delft, The Netherlands\\
E-mail: p.groeneboom@tudelft.nl
\end{document}